\documentstyle[epsf,amssymb,12pt]{amsart}
\input xy
\xyoption{all}

\newcommand{\la}{\langle}
\newcommand{\ra}{\rangle}

\newtheorem{theorem}{\bf Theorem}[section]
\newtheorem{lemma}[theorem]{\bf Lemma}
\newtheorem{lemmadef}[theorem]{\bf Lemma/Definition}
\newtheorem{remark}[theorem]{\bf Remark}

\newtheorem{corollary}[theorem]{\bf Corollary}
\newtheorem{definition}[theorem]{\bf Definition}

\DeclareFontFamily{OT1}{rsfs}{}
\DeclareFontShape{OT1}{rsfs}{n}{it}{<->rsfs10}{}
\DeclareMathAlphabet{\curly}{OT1}{rsfs}{n}{it}

\newcommand{\CC}{{\Bbb C}}

\newcommand{\NN}{{\Bbb N}}
\newcommand{\PP}{{\Bbb P}}

\newcommand{\RR}{{\Bbb R}}

\newcommand{\ZZ}{{\Bbb Z}}

\newcommand{\tlie}{{\frak t}}
\newcommand{\ulie}{{\frak u}}

\newcommand{\fX}{{\frak X}}

\newcommand{\IND}{\cong}
\newcommand{\DEND}{\sim}
\newcommand{\DDEND}{\simeq}

\newcommand{\Bihol}{\operatorname{Bihol}}

\newcommand{\ch}{\operatorname{ch}}

\newcommand{\End}{\operatorname{End}}
\newcommand{\EC}{E_{\CC}}

\newcommand{\Gr}{\operatorname{Gr}}
\newcommand{\Hom}{\operatorname{Hom}}
\newcommand{\Id}{\operatorname{Id}}

\newcommand{\Ker}{\operatorname{Ker}}
\newcommand{\Lie}{\operatorname{Lie}}

\newcommand{\Mor}{\operatorname{Mor}}
\newcommand{\parslope}{\operatorname{par-slope}}

\newcommand{\rk}{\operatorname{rk}}
\newcommand{\SL}{\operatorname{SL}}
\newcommand{\slope}{\operatorname{slope}}
\newcommand{\Spec}{\operatorname{Spec}}

\newcommand{\Tr}{\operatorname{Tr}}

\newcommand{\Vol}{\operatorname{Vol}}

\renewcommand{\Im}{\operatorname{Im}}

\newcommand{\ov}{\overline}
\newcommand{\un}{\underline}

\newcommand{\CCC}{{\cal C}}

\newcommand{\III}{{\curly I}}

\newcommand{\LLL}{{\curly L}}
\newcommand{\MMM}{{\curly M}}
\newcommand{\OOO}{{\curly O}}
\newcommand{\PPP}{{\curly P}}

\newcommand{\VVV}{{\cal V}}
\newcommand{\WWW}{{\cal W}}
\newcommand{\qu}{/\kern-.7ex/}
\newcommand{\exh}{\to\kern-1.8ex\to}

\newcommand{\uC}{\un{\CC}}
\newcommand{\ur}{\un{r}}
\newcommand{\ul}{\un{\lambda}}

\newcommand{\Gm}{\CC^{\times}}

\newcommand{\VP}{{\curly V}\kern-0.9ex\PPP}
\newcommand{\Veq}{{\curly V}_{\Gm}}
\newcommand{\VeqG}{{\curly V}_{\Gamma}}
\newcommand{\Veqr}{{\curly V}_{G(r)}}
\newcommand{\Veqi}{{\curly V}_{G(r_i)}}

\newcommand{\tV}{\widetilde{V}}

\newcommand{\imag}{\sqrt{-1}}

\author{Ignasi Mundet i Riera}
\address{Departamento de Matem{\'a}ticas, Universidad Aut{\'o}noma de Madrid,
Madrid, Spain}
\email{ignasi.mundet@@uam.es}
\date{January 31, 2001}

\title{Parabolic vector bundles and equivariant vector bundles}

\begin{document}
\maketitle
\begin{abstract}
Given a complex manifold $X$, a normal crossing divisor $D\subset X$ whose
irreducible components $D_1,\dots,D_s$ are smooth, and a choice
of natural numbers $\ur=(r_1,\dots,r_s)$, we construct a manifold
$X(D,\ur)$ with an action of a torus $\Gamma$ and we prove that some
full subcategory of the category of $\Gamma$-equivariant vector
bundles on $X(D,\ur)$ is equivalent to the category of parabolic
vector bundles on $(X,D)$ in which the lengths of the filtrations
over each irreducible component of $X$ are given by $\ur$. 
When $X$ is Kaehler, we study the Kaehler cone of $X(D,\ur)$ and
the relation between the corresponding notions of slope-stability.
\end{abstract}

\section{Introduction}

\subsection{}
Let $X$ be a (not necessarily compact) complex manifold and let 
$D\subset X$ be a divisor
with normal crossings, whose irreducible components $D_1,\dots,D_s$ are
smooth. Let $\PPP(X,D,\ur)$ be the category of parabolic vector bundles
over $(X,D)$ in which the lengths of the filtrations over the irreducible
components $D_1,\dots,D_s$ of $D$ are given by the natural numbers
$\ur=(r_1,\dots,r_s)$ (see Section \ref{ncd} for a precise definition).

In this paper we construct a manifold $X(D,\ur)$ endowed with an action of 
an algebraic torus $\Gamma$ and an invariant projection
$$\Pi:X(D,\ur)\to X,$$
we define a full subcategory $\VeqG(X(D,\ur))$
of the category of equivariant vector bundles over $X(D,\ur)$,
and we construct a functor 
$$M:\VeqG(X(D,\ur))\to\PPP(X,D,\ur),$$  
which induces an equivalence of categories (see Theorem \ref{mainthmur}).

When $X$ is Kaehler we describe the Kaehler cone of $X(D,\ur)$ in 
terms of that of $X$. Afterwards, for any choice of Kaehler class 
$\omega$ of $X$ and parabolic weights $\Lambda$ we construct a family 
of Kaehler classes
$\Omega(\omega,\Lambda,\epsilon)$ on $X(D,\ur)$ parametrized by
$\epsilon\in\RR_{>0}$ and we prove that for small enough $\epsilon$
the notions of $(\omega,\Lambda)$-parabolic slope stability and that
of $\Omega(\omega,\Lambda,\epsilon)$-slope stability for equivariant 
vector bundles correspond each other by the equivalence of categories
(see Theorem \ref{estabilitat}).
For any set of weights $\Lambda$ we have 
$$\lim_{\epsilon\to 0}\Omega(\omega,\Lambda,\epsilon)=\Pi^*\omega$$
so when $\epsilon\to 0$ the volume of the fibres of $\Pi$ tends to
$0$. In this sense, our result could be thought of as a statement
on stability of vector bundles in the adiabatic limit.

The functor $M$ gives also an equivalence of categories of families
of (parabolic, equivariant) vector bundles parametrized by complex
spaces $S$, and in some particular cases our results allow to identify the two
moduli problems.

The construction of $X(D,\ur)\to X$ is functorial
with respect to $X$ and $D$ in the following sense: if $f:Y\to X$ is a 
map of complex manifolds which is transverse to $D$, then we have
an induced map $f_{D,\ur}$ which makes the following diagram commutative:
$$\xymatrix{Y(f^{-1}D,\ur)\ar[r]^-{f_{D,\ur}}\ar[d] & X(D,\ur)\ar[d] \\
Y \ar[r]^{f} & X.}$$
Furthermore the functor $M$ is compatible with the maps
$f_{D,\ur}$, in the sense that the following
diagram is also commutative:
$$\xymatrix{\VeqG(X(D,\ur)) \ar[r]^-{f_{D,\ur}^*}\ar[d]_M & 
\VeqG(Y(f^{-1}D,\ur)) \ar[d]^M \\
\PPP(X,D,\ur) \ar[r]^{f^*} & \PPP(Y,f^{-1}D,\ur).}$$

The category $\VeqG(X(D,\ur))$ consists of equivariant vector bundles
$W\to X(D,\ur)$ whose weights over the fixed point set of the action 
of $\Gamma$ on $X(D,\ur)$ satisfy certain restrictions. Hence, the
condition of some equivariant vector bundle $W$ being an object of
$\VeqG(X(D,\ur))$ is purely topological.

\subsection{}
The ideas in this paper are very similar of those of Biswas in
\cite{B} (see also the references therein). 
Biswas constructs a finite Galois covering $Y$ of $X$ and shows 
how to obtain, out of a parabolic bundle over $(X,D)$,
a bundle on $Y$ which is equivariant w.r.t. the Galois group of
the covering (this establishes a link between the theory of parabolic 
bundles and the theory of bundles over orbifolds --- note that some
particular cases of this link were already known before the work of
Biswas). Biswas also relates the stability
conditions of the parabolic sheaf and of the equivariant one. 
In contrast with our case, however, the manifold $Y$ depends not only on 
$X$, $D$ and $\ur$, but also on the choice of the parabolic weights. 
It would be interesting to relate the approach of Biswas to that of this paper.

From another point of view, the results which we present here are 
related to those of Garc\'{\i}a--Prada \cite{GP1,GP2}, in which holomorphic
pairs are studied in terms of equivariant vector bundles. The ideas of
Garc\'{\i}a--Prada have been successfully applied to other situations 
dealing with holomorphic bundles with extra structure (see
for example \cite{A, AGP,BDGW, BGM}).

\subsection{}
The construction of $X(D,\ur)$ is made in several steps. First we consider
the case of $D$ smooth and $r=r_1=1$. We define $X(D,1)$ (or $X_D$
for short) as a family of conics over $X$ which degenerate precisely over
$D$. This is the same thing as the blow up of $X\times\PP^1$ along 
$D\times\{[0:1]\}$, which is (a compactification of) the deformation to the normal
cone of $D$. 
The manifold $X_D$ inherits an action from the diagonal action of $\CC^*$
on $X\times\PP^1$ which is trivial on the $X$ factor and which is defined
on $\PP^1$ as $\theta\cdot[y:w]:=[y:\theta w]$. 
Next we consider the case of $D$ smooth and $r>1$. 
We construct a tower of manifolds 
$$X(D,r)=Y_r\to Y_{r-1}\to\dots\to Y_0=X$$
by applying recursively the previous construction, so that 
$Y_{j+1}=(Y_j)_{D_j}$. Here $D_j\subset Y_j$ is a certain smooth divisor which lies
over $D$, and each manifold $Y_j$ carries an action of $(\Gm)^j$.

Finally, when $D=D_1\cup\dots\cup D_s$ is a normal crossing divisor, we define
$X(D,\ur)$ to be the fibred product
$$\Pi:X(D,\ur)=X(D_1,r_1)\times_X\dots\times_X X(D_s,r_s)\to X.$$
This carries a diagonal action of $\Gamma=G(r_1)\times\dots\times G(r_s)$.
The fibres of $\Pi$ over $X\setminus D$ are products of $\PP^1$'s,
and those over $D$ are singular, their irreducible components
being products of $\PP^1$'s. In fact, $X(D,\ur)$ can be constructed
by making a sequence of blow ups along subvarieties of a product
$X\times\PP^1\times\dots\PP^1$ which lie above $D$.

The category $\VeqG(X(D,\ur))$ and the
functor $M$ are also defined inductively, following the construction of
$X(D,\ur)$. The main point is to define $M$ in the case $D$ smooth and
$r=1$, since everything is built up by applying recursively the construction in
this simplest case. The proof that $M$ induces an equivalence of categories and
some of the computations needed to prove the relation with Mumford--Takemoto 
stability are also obtained by reducing to this simple case.

\subsection{}
We now explain the contents of the following sections.
In Section \ref{sec:equibdls} we give some definitions and results
on equivariant vector bundles which will be used along the paper.
In Section \ref{construccioXD} we describe the construction of 
the $\Gm$-manifold $X_D$ out of the pair $(X,D)$; as we said before,
we will obtain the manifold $X(D,\ur)$ by iterating this construction.
In Section \ref{sec:equicat} we define a functor from (some full
subcategory of) the category of $\Gm$-bundles over $X_D$ to that
of parabolic vector bundles over $(X,D)$ with $r=r_1=1$, and we prove
that it induces an equivalence of categories.
In Section \ref{divllis} we extend the result of the previous section
to arbitrary parabolic bundles over a smooth divisor.
Section \ref{sec:cohoques} is devoted to studying the Kaehler cone
of the manifolds $X(D,r)$ and some topological aspects of the equivalence
between equivariant and parabolic vector bundles.
These results are used in Section \ref{slopesmooth} to relate the notions
of slope and parabolic slope of vector bundles.
In Section \ref{ncd} everything is extended to the case of a normal
crossing divisor. 
In Section \ref{sec:stab} we study the stability condition.
Finally, in Section \ref{sec:dim1} we consider the case of $X$
being a Riemann surface.

\subsection{Notations and conventions.}
Unless we say the contrary, the following will be implicitly assumed in
this paper: all vector bundles will be complex, all metrics on vector
bundles will be Hermitian, all vector bundles, maps of vector bundles,
manifolds, and actions of groups on manifolds will be holomorphic.
A divisor will mean a reduced divisor.
 
We will use the following notations.
If $Y$ is a complex manifold $\Bihol(Y)$ will denote the group of
biholomorphisms of $Y$ with itself.
If $G$ is an abelian group, $Y$ is a $G$ manifold, $V\to Y$ is a $G$-equivariant
vector bundle and $Y'\subset Y^G$, then $\chi_G(V|_{Y'})$ 
will denote the set of characters appearing in the decomposition of
$V|_{Y'}$ in bundles of irreps of $G$ (as usual the superscript $^G$
denotes the fixed points). We will call $\chi_G(V|_{Y'})$
the set of $G$-weights of $V$ on $Y'$.

If $\CCC$ is any category, we will usually write $A\in\CCC$ to mean
that $A$ is an object of $\CCC$.

\subsection {Aknowledgements.}
Most of this paper was done during a stay at the Centre de Math\'ematiques
of the \'Ecole Polytechnique. The author is very pleased to have the opportunity 
to thank this institution for financial support and for providing excellent 
working conditions. Thanks are also due to C. Sabbah, V. Navarro Aznar for 
some useful comments.

\section{Equivariant bundles}
\label{sec:equibdls}

\subsection{Weights of bundles}
\label{def:stronginclusion}
We will say that an inclusion $Y\subset Z$ of topological spaces is a
strong inclusion if every connected component of $Z$ contains some
point of $Y$.
Assume that a group $G$ acts on a vector bundle $W\to Z$ linearly on the 
fibres. If $Z'\subset Z^G$, and $Y\subset Z'$ is a strong inclusion, 
then $$\chi_G(W|_Y)=\chi_G(W|_Z).$$ This follows from the local
invariance of the set of weights on the fixed point locus.

\subsection{Stable and unstable sets}
\label{Gmbundles}

Let $Y$ be a manifold with an action of $\Gm$. Denote the image 
of $\theta\in\Gm$ acting on $y\in Y$ by $\theta\cdot y$.

Let $Y'$ be a connected component of $Y^{\Gm}$. Define the stable
(resp. unstable) set $U^+(Y')$ (resp. $U^-(Y')$) of $Y'$ to be
\begin{align*}
U^+(Y') &=\{y\in Y\mid \lim_{\Gm\ni\theta\to 0}\theta\cdot y\text{
  exists and belongs to $Y'$}\}, \\
U^-(Y') &=\{y\in Y\mid \lim_{\Gm\ni\theta\to\infty}\theta\cdot y\text{
  exists and belongs to $Y'$}\}.
\end{align*}
One can prove that the sets $U^+(Y')$ and $U^-(Y')$ are smooth complex
submanifolds of $Y$ (this is part of the Kaehler version of
Bia\l ynicki-Birula theorem, see \cite{CS}). In
the cases which will be considered in this paper this fact will be
obvious. 

Let $W\to Y$ be a complex bundle with a linear action of $\Gm$
lifting the one on $Y$. We will denote the fibre over $y$ by $W_y$.
For any fixed point $y\in Y^{\Gm}$ the fibre $W_y$ gets a linear
action of $\Gm$. The isomorphism class of this 
representation depends only on the connected
component $Y'$ of $Y^{\Gm}$ to which $y$ belongs. Let
$W_y=\bigoplus_{k\in\ZZ}W_y(k)$ be the decomposition of $W_y$ 
in weights. 

\begin{lemma}
Consider the standard inclusion of groups $S^1\subset \Gm$, and take
on $Y$ the action of $S^1$ induced by restriction.
Let $y\in Y^{S^1}$ be a fixed point. There exists a $S^1$ invariant
neighbourhood $A\subset Y$ of $y$ and a holomorphic and
$S^1$-equivariant trivialisation
$\phi:W|_A\stackrel{\simeq}{\longrightarrow} 
A\times W_y$, with $A\times W_y$ supporting the diagonal
action of $S^1$. (We will say that such a trivialisation is
centered around $y$.)
\label{modellocalW}
\end{lemma}
\begin{pf}
Let $A_0$ be a $S^1$-invariant
neighbourhood of $y$, small enough so that there is a
holomorphic trivialisation of $W|_{A_0}$. Then the map
$e:\Gamma:=\Gamma(A_0;W)\to W_y$ given by evaluation at $y$ (where
$\Gamma(A_0;W)$ denotes the set of holomorphic sections defined on
$A_0$) is exhaustive. The group $S^1$ acts on $\Gamma$ by
pullback. Let $\Gamma=\bigoplus_{k\in\ZZ}\Gamma(k)$ be the splitting
given by the weights. Since $e$ is exhaustive,
there exists elements $w_1,\dots,w_N$ (where
$N=\rk W$) and weights $k_1,\dots,k_N$
such that $w_j\in\Gamma(k_j)$ and such that
$e(w_1),\dots,e(w_N)$ span $W_y$ (and so form a basis). Consequently,
there is a neighbourhood $A\subset A_0$ of $y$ over which 
the sections $\{w_j\}$ trivialise $W$.
\end{pf}
\begin{remark}
Note that $Y^{S^1}=Y^{\Gm}$. The obtained trivialisation is 
weakly $\Gm$-equivariant in the following sense. 
If $y\in A$, $\theta\in\Gm$ and $\theta$ can be joined to $1\in
\Gm$ by a path $\{\theta(t), 0\leq t\leq 1\}$ so that 
$\theta(t)\cdot y\in A$ for all $0\leq t\leq 1$, then 
for any $v\in W_y$ we have $\phi(\theta\cdot
v)=\theta\cdot\phi(v)$. This follows from the fact that both the 
trivialisation $\phi$ and the action of $\Gm$ on $W$ are holomorphic.
\label{modellocalWrem}
\end{remark}

Take a metric on $W$. 
Chose a connected component $Y'$ of the fixed point set and
define, for any $r\in\RR$, 
$$W^{+,r}(Y') = \{w\in W|_{U^+(Y')}\mid \lim_{\theta\to 0}
|\theta^{-r}(\theta\cdot w)|<\infty\}.$$
We call $W^{+,r}(Y')$ the $r$-stable subbundle of $W$ towards $Y'$.

\begin{lemma}
The family of sets  $\{W^{+,r}(Y')\mid r\in\RR\}$ 
is a decreasing filtration of $\Gm$-invariant subbundles
of $W|_{U^+(Y')}$, and it is independent of the chosen metric on $W$.
Let $\phi:W\to W'$ be an equivariant map of $\Gm$-bundles. Then, for
any $r\in\RR$, we have $\phi(W^{+,r}|_{U^+(Y')})\subset
{W'}^{+,r}|_{U^+(Y')}$.
\label{inclusio}
\end{lemma}
\begin{pf}
Follows from the definitions and Lemma \ref{modellocalW} together with
Remark \ref{modellocalWrem}.
\end{pf}


Suppose that the only weight of $W$ in $Y'$ is zero. 
For any $z\in U^-(Y')$ define the map $\rho^z_W:W_z\to W_{z_-}$ by
$\rho^z_W(w):=\lim_{\theta\to\infty}\theta\cdot w$.
Similarly, if $z\in U^+(Y')$, we set $R^z_W:W_z\to W_{z_+}$ to be
$R^z_W(w):=\lim_{\theta\to 0}\theta\cdot w$.

\begin{lemma}
If $z\in U^-(Y')$ (resp. $z\in U^+(Y')$) then $\rho^z_W$
(resp. $R^z_W$) is well defined and is an isomorphism of vector spaces.
\label{rhoiso}
\end{lemma}
\begin{pf}
This also follows from Lemma \ref{modellocalW} and Remark
\ref{modellocalWrem}. 
\end{pf}

\subsection{Closure of subbundles}
\label{subs:closure}
Let $Z$ be a manifold, let $Y=Z\times\CC$, and
consider the action of $\Gm$ on $Y$ given by
$\theta\cdot(z,a)=(z,\theta a)$. The 
fixed point set of this action is $Y'=Z\times\{0\}$, and we have
$U^+(Y')=Y$. Define $Y^*=Y\setminus Y'$. Let $p:Y\to Z$ be the
projection. 

\begin{lemma}
Let $W\to Y$ be a $\Gm$-vector bundle whose only weight in $Y'$ is
$0$. Let $W'\subset W|_{Y^*}$ be a $\Gm$-invariant
subbundle. Then there is a unique $\Gm$-invariant
subbundle $\ov{W'}\subset W$ such that $\ov{W'}|_{Y^*}=W'$.
Furthermore, if $V\to Y$ is another $\Gm$-vector bundle whose only
weight at $Y'$ is $0$ and $\phi:W\to V$ is a $\Gm$-equivariant map
such that $\phi|_{W'}:W'\to V|_{Y^*}$ is exhaustive, then 
$\phi|_{\ov{W'}}:\ov{W'}\to V$ is also exhaustive.
\label{closuresbdl}
\end{lemma}
\begin{pf}
We first prove uniqueness. Let $\ov{W'}\subset W$ be such an
extension. Necessarily, the only weight of $\ov{W'}$ in $Y'$ is $0$.
Let $z_0\in Z$ be any point and define $z=(z_0,1)\in Y$, so
that $z_+=(z_0,0)$. By Lemma \ref{rhoiso}, $R^z_{\ov{W'}}$ 
is an isomorphism and by
equivariance we have $R^z_{\ov{W}'}=R^z_W|_{W'_z}$. Hence we must have
$\ov{W'}_{(z_0,0)}=R^z_W(W'_z)\subset W_{(z_0,0)}$. In other words, the
fibre of $\ov{W'}$ over $(z_0,0)$ is determined by $R^z_W$ and by the
fibre $W'_z$. The last claim of the lemma follows also easily using
the isomorphisms $R^z$.

To prove existence it suffices to work locally (thanks to uniqueness).
So take any $z_0\in Z$. By Lemma \ref{modellocalW} there exists a
neighbourhood $A\subset Y$ of $(z_0,0)$ and a weakly $\Gm$-equivariant
trivialisation 
$$\phi:W|_A\to A\times W_{(z_0,0)},$$ where $\Gm$ acts trivially on
$W_{(z_0,0)}$. By shrinking $A$ if necessary, we may assume that 
$$A=\{(b,a)\in Z\times\CC \mid b\in p(S),\ |a|<\epsilon\}$$ for some
$\epsilon>0$. Let $A^*=A\cap Y^*$.
Using the trivialisation $\phi$ the subbundle $W'|_{S^*}$ is described 
by a map to the Grassmannian of the fibre over $(z_0,0)$:  
$$\Psi:A^*\to \Gr(W_{(z_0,0)}).$$
Since $W'$ is $\Gm$-invariant and holomorphic, we deduce that $\Psi$
is weakly $\Gm$-equivariant and holomorphic, i.e., 
$\Psi(z,a)=\Psi(z,b)$ whenever $0<|a|<\epsilon$ and
$0<|b|<\epsilon$. It then follows that we can extend $\Psi$ to a map
$\ov{\Psi}:A\to \Gr(W_{(z_0,0)})$ by simply setting
$\ov{\Psi}(z,a)=\Psi(z,\epsilon/2)$. This map describes the desired
weakly $\Gm$-equivariant subbundle of $W|_A$. 
\end{pf}

The preceding lemma holds also true when $W$ has a unique weight on
$Y'$ (not necessarily zero). However, if $W$ has more than one
weight in $Y'$ then the result might be false. Indeed, consider the
trivial bundle $\un{\CC^2}\to\CC^2$ (i.e., here $Z=\CC$). 
Let $\la e_1,e_2\ra$ be the
canonical basis of the fibre $\CC^2$. Let $L\subset\un{\CC^2}$ be the
line bundle defined over $\CC\times\Gm$ whose fibre over
$(z,a)$ is $L_{(z,a)}=\la z e_1+a e_2\ra$. Consider the action of $\Gm$
on $\CC^2$ given by $\theta\cdot(z,a)=(z,\theta a)$ and its lift to
the trivial bundle $\un{\CC^2}$ defined by $\theta\cdot(\alpha
e_1+\beta e_2)=\alpha e_1+\theta \beta e_2$. Then $L$ is 
$\Gm$-invariant, but it does not extend to a line subbundle of
$\un{\CC^2}\to\CC^2$. 

\begin{corollary}
Let $W_i\to Y$, $i=0,1$, be $\Gm$-equivariant bundles whose only
weight in $Y'$ is $0$, and let
$\phi:W_0|_{Y^*}\to W_1|_{Y^*}$ be a $\Gm$-equivariant map. 
Then $\phi$ extends to a $\Gm$-equivariant map from $W_0$ to $W_1$.
\label{corclosuresbdl}
\end{corollary}
\begin{pf}
Let $W=W_0\oplus W_1$ and let $W'$ be the graph of
$\phi$. The only weight of $W$ in $Y'$ is $0$, and
$W'$ is a $\Gm$-equivariant subbundle of $W|_{Y^*}$. Hence by Lemma
\ref{closuresbdl} $W'$ extends to a $\Gm$-equivariant
subbundle $\ov{W'}\subset W$. To see
that $\ov{W'}$ is the graph of a map from $W_0$ to $W_1$ it suffices to
check that the projection $\ov{W'}\to W_0$ is an isomorphism,
and this follows from the last claim of Lemma \ref{closuresbdl}.
\end{pf}

\subsection{Lifting actions to line bundles}

Let $X$ be a manifold, and let $D\subset X$ be a smooth
divisor. Let $\pi:L\to X$ be a line bundle and $\sigma\in H^0(X;L)$ a
section which is transverse to the zero section and such that
$\sigma^{-1}(0)=D$. 

\begin{lemmadef}
Let $\Bihol(X,D)$ be the group of biholomorphisms of $X$ which
preserve $D$, and let $\Bihol(L)$ be the group of biholomorphisms of
$L$ which preserve $\pi^{-1}(D)$ and which map fibres linearly to
fibres. There exists a right inverse
$\phi_{\sigma}:\Bihol(X,D)\to\Bihol(L)$ to the projection map
$\Bihol(L)\to\Bihol(X,D)$.
\label{useful1}
\end{lemmadef}
\begin{pf}
Let $f\in\Bihol(X,D)$. We first define the restriction $f_0$ of
$\phi_{\sigma}(f)$ to $L_0:=L\setminus\pi^{-1}(D)$ as follows: for any
$x\in L_0$ we set
$$f_0(x)=\frac{x}{\sigma(\pi(x))}\sigma(f\pi(x))\in L_{f\pi(x)}$$
(note that the fraction in the RHS is a complex number). The map
$f_0$ can be equivalently seen as a map of line bundles
$$F:L|_{X\setminus D}\to f^*L|_{X\setminus D},$$
and to see that $f_0$ extends to a map from $L$ to $L$ it
suffices to prove that $F$ extends to a map of line bundles
defined over the whole $X$. 
For that, and thanks to Riemann's extension theorem, it is enough to
check that if we fix a metric on $L$ and
$K\subset X$ is a compact subset, then the restriction of $F$ to
$K\cap(X\setminus D)$ is bounded. This is the same thing as 
$$\sup\left\{\frac{|\sigma(fy)|}{|\sigma(y)|}\ :\  y\in
  K\cap(X\setminus D)\right\} < \infty,$$ and this follows from the
fact that $f$ is a biholomorphism which preserves $D$ and that
$\sigma$ is transverse to the zero section. Hence $f_0$ extends to a
map $\phi_{\sigma}(f):L\to L$ which lifts $f$.
It is also clear that $\phi_{\sigma}(\Id)=\Id$, and that if
$f,g\in\Bihol(X,D)$ then $\phi(s)(fg)=\phi_{\sigma}(f)\phi_{\sigma}(g)$
(indeed, that $f_0g_0=(fg)_0$ is obvious from the definition, and then
use that $X\setminus D$ is dense in $X$). Finally, it follows from
the lattest property that if $f\in\Bihol(X,D)$ then 
$\phi_{\sigma}(f)\phi_{\sigma}(f^{-1})=\Id$,
so $\phi_{\sigma}(f)$ is indeed a biholomorphism. (Note that
$\phi_{\sigma}(f)$ can be characterized as the unique element of
$\Bihol(L)$ which lifts $f$ and which makes $\sigma$ equivariant.) 
\end{pf}

\begin{lemma}
If $f\in\Bihol(X,D)$, then the restriction of $\phi_{\sigma}(f)$ to
$\pi^{-1}(D)$ can be computed as follows. Let $z\in D$. Since $\sigma$
is transverse to zero, $d\sigma(z):N_z\to L_z$ is an isomorphism,
where $N\to D$ is the normal bundle of $D\subset X$. On the other
hand, $df:N_z\to N_{f(z)}$ is also an isomorphism and, furthermore,
that $$\phi_{\sigma}(f)|_{L_z}=d\sigma(f(z))\circ df\circ
d\sigma(z):L_z\to L_{f(z)}.$$
\label{useful3}
\end{lemma}
\begin{pf}
It follows from a computation in local coordinates.
\end{pf}
 
\begin{lemma}
If a group $\Gamma$ acts on $X$ by biholomorphisms preserving $D$,
then there is a unique linear lift of the action of $\Gamma$ to $L$
for which $\sigma$ is $\Gamma$-equivariant. If $x\in X^{\Gamma}$, the
character $\chi$ of $\Gamma$ acting on $L_x$ is $1$ if $\sigma(x)\neq
0$, and if 
$\sigma(x)=0$ then $\chi$ is equal to the character of the action of $\Gamma$
on the fibre $N_x$ of the normal bundle $N\to\sigma^{-1}(0)$.
\label{useful2}
\end{lemma}
\begin{pf}
The first part of the lemma follows from Lemma \ref{useful1}, and the
second part from Lemma \ref{useful3}.
\end{pf}

\section{The $\Gm$-manifold $X_D$}

\label{construccioXD}
Let $X$ be a smooth manifold (not necessarily compact), and let
$D\subset X$ be a smooth divisor.
Let $L\to X$ be a line bundle and $\sigma\in H^0(L)$ a section
which is transverse to zero and such that $\sigma^{-1}(0)=D$.
Let 
$$X_{D,\sigma}:=\{[x:y:w]\in\PP(L\oplus\uC\oplus\uC)_z\mid z\in X,\
xy=w^2\sigma(z)\},$$ 
and define an action of $\Gm$ on $X_D$ by
$\theta\cdot [x:y:w] := [\theta^2 x:y:\theta w]$
for any $\theta\in\Gm$.
If $\sigma'\in H^0(L)$ is another nonzero section which is transverse
to zero then $\sigma'=\theta\sigma$ for some $\theta:X\to\Gm$, so the
map $\Theta:X_{D,\sigma}\to X_{D,\sigma'}$ which sends $[x:y:w]$ to
$[\theta x:y:w]$ is a biholomorphism. In view of this, we will just
write, to save on notation, $X_D$ instead of $X_{D,\sigma}$. 
(But it is important to keep in mind that whereas the assignement 
$(X,D,\sigma)\to X_{D,\sigma}$ is functorial, this is not the case 
of the assignement $(X,D)\to X_D$.) 
We will denote by $p:X_D\to X$ the projection.

It follows from an easy local computation that $X_D$ is smooth.

\begin{figure}[htb]
\epsfxsize=10cm
\epsffile{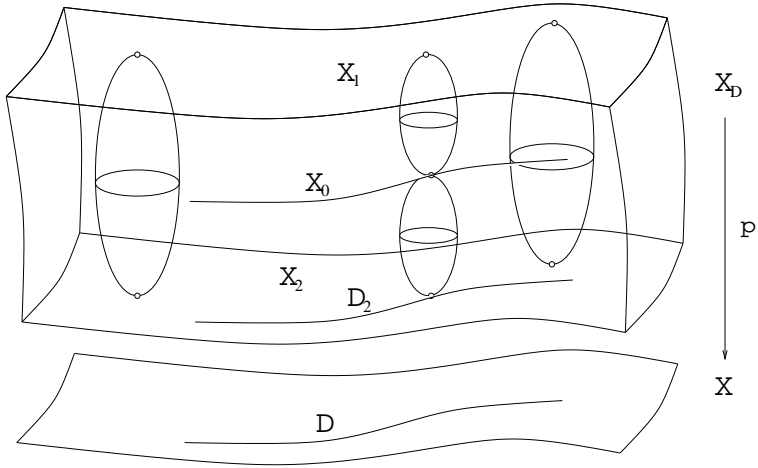}
\caption{The manifold $X_D$. The spheres represent the closure of
orbits of the action of $\Gm$.}
\end{figure}

\begin{lemmadef}
There exists a morphism of groups $\phi$ from $\Bihol(X,D)$ to
$\Bihol(X_D)$ such that for any $f\in\Bihol(X,D)$ we have
$p\circ\phi(f)=f\circ p$.
\label{useful}
\end{lemmadef}
\begin{pf}
Use Lemma \ref{useful1} to get a lift $\phi_{\sigma}(f):L\to L$ of
$f$. Taking the trivial lift of the action to the trivial bundle
$\uC$, this gives an element of $\Bihol(\PP(L\oplus\uC\oplus\uC))$ which
leaves $X_D$ fixed, and we define $\phi(f)$ to be the induced element
of $\Bihol(X_D)$. 
\end{pf}

In the sequel we will denote, for any $f\in\Bihol(X,D)$, $f_D:=\phi(f)$. 

The fixed point set of the action of $\Gm$ on $X_D$ is the disjoint
union of the following submanifolds:
$$X_1 = \{[x:y:w]=[0:1:0]\},\qquad
X_2 = \{[x:y:w]=[1:0:0]\},$$
$$X_0 = \{[x:y:w]=[0:0:1]\}\cap X_D.$$
There are canonical identifications given by the projection $p$:
$$X_2\simeq X\simeq X_1\qquad\text{ and }\qquad X_0\simeq D.$$ Using
this identifications we have, for any $f\in\Bihol(X,D)$, 
$$f_D|_{X_1}=f_D|_{X_2}=f\qquad\text{ and }\qquad
f_D|_{X_0}=f|_D.$$
Denote by $D_2\subset X_D$ the copy of $D$ obtained by
means of the identification $X_2\simeq X$.

Define $\Delta^+=\Delta^+(X,D):=\ov{U^+(X_0)}$ and 
$\Delta^-=\Delta^-(X,D):=\ov{U^-(X_0)}$.
Both $\Delta^+$ and $\Delta^-$ are 
smooth $\Gm$-invariant divisors in $X_D$, and  
$\Delta^+\cap X_2=D_2$.
We can explicitly describe them as
\begin{equation}
\Delta^+=\{y=0\}\cap X_D
\qquad\text{ and }\qquad
\Delta^-=\{x=0\}\cap X_D.
\label{descriuD}
\end{equation}
Let us denote by $N^+\to\Delta^+$ (resp. $N^-\to\Delta^-$) the normal
bundle of the inclusion $\Delta^+\subset X_D$ (resp. $\Delta^-\subset
X_D$).

\begin{lemma}
The fixed point set $(\Delta^+)^{\Gm}$ is the disjoint union of
$X_0=\Delta^+\cap X_0$ and $X_2':=\Delta^+\cap X_2$. The weight of the
action of $\Gm$ on the restriction of $N^+$ to $X_0$ (resp. $X_2'$) is
$-1$ (resp. $0$).
\label{weightN}
\end{lemma}
\begin{pf}
The first statement is obvious.
To prove the statement on the weights, observe that since $\Delta^+$ and
$X_2$ intersect transversely, $N^+|_{X_2'}$ can be identified (in a
$\Gm$-equivariant way) with the normal bundle of the inclusion
$X_2'\subset X_2$. But since the action of $\Gm$ on $X_2$ is trivial,
it follows that the weight of $\Gm$ acting on $N^+|_{X_2'}$ is $0$. 
On the other hand, $N^+|_{X_0}$ can be identified with the normal bundle
of the inclusion $X_0\subset U^-(X_0)$. Now, since
$U^-(X_0)=\{[0:y:1]\}\cap X_D$, and
$\theta\cdot[0:y:1]=[0:y:\theta]=[0:\theta^{-1}y:1]$, it follows that
the action of $\Gm$ has weight $-1$.
\end{pf} 

\begin{lemma}
The fixed point set $(\Delta^-)^{\Gm}$ is the disjoint union of
$X_0=\Delta^-\cap X_0$ and $X_1':=\Delta^+\cap X_1$. The weight of the
action of $\Gm$ on the restriction of $N^-$ to $X_0$ (resp. $X_1'$) is
$1$ (resp. $0$).
\label{weightN2}
\end{lemma}
\begin{pf}
Exactly like that of the preceding lemma, but using
that $N^-|_{X_0}$ can be identified with the normal bundle
of the inclusion $X_0\subset U^+(X_0)$, and that
$U^+(X_0)=\{[x:0:1]\}\cap X_D$ and
$\theta\cdot[x:0:1]=[\theta^2x:0:\theta]=[\theta x:0:1]$.
\end{pf}

\subsection{Example: The manifold $H_n$}
\label{defHn}
Define for any natural number $n\geq 1$
$$H_n=\{(t_1,\dots,t_n,[x:y:w])\mid
xy=w^2t_1\}\subset\CC^n\times\PP(\CC^3).$$ 
Consider the action of $\Gm$ on $H_n$ given by
$\theta\cdot(t_1,\dots,t_n,[x:y:w])=
(t_1,\dots,t_n,[\theta^2x:y:\theta w]).$  
Let $p_n:H_n\to\CC^n$ be the projection. 
Note that the $\Gm$-manifold $H_n$ is nothing but 
$\CC^n_{\{0\}\times\CC^{n-1}}$. 

For future use we define the map $\xi:H_n\to\RR$ as follows
$$\xi(t_1,\dots,t_n,[x:y:w]):=
-\frac{2|x|^2+|w|^2}{|x|^2+|y|^2+|w|^2}.$$ 
Up to some multiplicative constant, this is the moment map of the
action of $S^1\subset\Gm$ on $H_n$ with respect to the symplectic
structure inherited by the inclusion $H_n\subset\CC^n\times\PP(\CC^3)$ and
the usual symplectic structures on $\CC^n$ and $\PP(\CC^3)$.
Finally, we define the following subset of $H_n$:
$$U_n = \{(t_1,\dots,t_n,[x:y:w])\in H_n\mid x\neq
0\}=U^+(\{[x:y:w]=[1:0:0]).$$  

\section{An equivalence of categories}
\label{sec:equicat}

\subsection{Definitions of $\Veq(X_D)$, $\PPP(X,D)$ and $\mu$}
We define the category $\Veq(X_D)$ as follows:
\begin{enumerate}
\item The objects of $\Veq(X_D)$ are $\Gm$-equivariant vector
  bundles $W\to X_D$ such that the only weight in $W|_{X_1\cup X_2}$ is $0$ 
  and the weights in $W|_{X_0}$ belong to $\{0,1\}$.
\item The morphisms between two $\Gm$-equivariant bundles $W,W'$ are
  the $\Gm$-equivariant maps of vector bundles $\psi:W\to W'$.
\end{enumerate}
We define the category $\PPP(X,D)$ as follows:
\begin{enumerate}
\item The objects of $\PPP(X,D)$ are pairs $(V,V_1)$ consisting of a
vector bundle $V\to X$ together with a subbundle $V_1\subset V|_D$.
\item The morphisms between two objects $(V,V_1)$ and $(V',V'_1)$
  are the maps of vector bundles $\phi:V\to V'$ such that
  $\phi|_D(V_1)\subset V'_1$.
\end{enumerate}

\subsubsection{The functor $\mu$}
\label{defmu}
We are now going to define a functor $\mu:\Veq(X_D)\to \PPP(X,D)$. 
By a slight abuse of notation, in this subsection we will
write $D=D_2$ and $X=X_2$. 

We first define $\mu$ acting on objects.
Let $W\to X_D$ be a $\Gm$-bundle. Let
$W^+:=W^{+,1}(X_0)\to U^+(X_0)$ be the $1$-stable subbundle of $W$
towards $X_0$ (see Subsection \ref{Gmbundles}). 
The bundle $W^+$ extends to a unique
$\Gm$-invariant subbundle $\ov{W^+}\subset W|_{\Delta^+}$. 
Indeed, locally around any $z\in D\subset X_D$, $X_D$ is $\Gm$-equivariantly 
biholomorphic to $Y=Z\times\CC$ with $Z=\CC^n$ and the action of 
$\Gm$ described in \ref{subs:closure}.
Hence, by Lemma \ref{closuresbdl}, there exists the closure
$\ov{W^+}\subset W$ of $W^+$ near $z$. 
By uniqueness we can make the same reasoning around all $z\in D$ and
patch the resulting local closures of $W^+$, thus getting the desired
extension $\ov{W^+}$. Now, we define 
$$\mu(W):=(W|_X,\ov{W^+}|_D).$$

We state the following lemma for later use.

\begin{lemma}
Let $W\in \Veq(X_D)$, and let $W_0\subset W|_{\Delta^+}$ be a $\Gm$-equivariant 
subbundle, and assume that $\chi_{\Gm}(W_0|_{X_0})=\{1\}$. Let $(V,V_1)=\mu(W)$. 
Then, identifying $D$
with $\Delta^+\cap X_2$ we have an inclusion of vector bundles over $D$:
$$W_0|_{\Delta^+\cap X_2}\subset V_1.$$
\label{inclsubbdl}
\end{lemma}
\begin{pf}
It follows easily from the definition of $\mu$.
\end{pf}
To define $\mu$ acting on morphisms, observe that if $\psi:W\to W'$ is
a $\Gm$-equivariant map of vector bundles then Lemma \ref{inclusio} implies
that $\psi(W^+)\subset {W'}^+$. From this it follows easily that 
$\psi(\ov{W^+})\subset \ov{{W'}^+}$. Hence, the restriction of $\psi$ to $X$ is
a morphism in the category $\PPP(X,D)$ between $\mu(W)$ and
$\mu(W')$, and we define
$$\mu(\psi):=\psi|_X\in\Mor_{\PPP(X,D)}(\mu(W),\mu(W')).$$ 

\subsubsection{}
\label{functorialitat}
It is straightforward to deduce from the previous definition that if 
$W,W',W''\in\Veq(X_D)$ and $\psi\in\Mor_{\Veq(X_D)}(W,W')$ and
$\xi\in\Mor_{\Veq(X_D)}(W',W'')$ then 
$$\mu(\xi\circ\psi)=\mu(\xi)\circ\mu(\psi).$$
This proves that $\mu$ is indeed a functor.

The following is easily checked: if $f:Y\to X$ is a map whose image is
transverse to $D$ (so that $f^{-1}(D)\subset Y$ is a smooth divisor)
and $f_D:Y_{f^{-1}(D)}\to X_D$ is the map induced by taking on $Y$
the pullback of $L$ and $\sigma\in H^0(L)$, then we get two
commutative diagrams:
\begin{equation}
\xymatrix{
Y_{f^{-1}D}\ar[r]^{f_D}\ar[d] & X_D\ar[d] \\
Y \ar[r]^{f} & X,}
\qquad\qquad
\xymatrix{
\Veq(X_D)\ar[d]_{\mu}\ar[r]^-{f_D^*} & \Veq(Y_{f^{-1}(D)})\ar[d]^{\mu} \\
\PPP(X,D)\ar[r]^-{f^*} & \PPP(Y,f^{-1}(D)),}
\label{commutamu}
\end{equation}
where the horizontal arrows in the RHS are the functors induced by pullbacks.
Note that a particular case of such a map $f:Y\to X$ is any
biholomorphism $f\in\Bihol(X)$. 

\begin{theorem} 
The functor $\mu:\Veq(X_D)\to\PPP(X,D)$ induces an equivalence of
categories. 
\label{mainthm}
\end{theorem}

\begin{remark}
It is natural to ask in view of the preceding theorem 
what happens if we consider $\Gm$-equivariant bundles on $X_D$ which
do not belong to $\Veq(X_D)$. A partial answer to this
question is the following: the category consisting
of $\Gm$-equivariant bundles whose weights on 
$X_D|_{X_1\cup X_2}$ are $0$ and whose weights in $X_D|_{X_0}$ belong
to $\{0,k\}$ for some $k\in\NN$ (instead of $\{0,1\}$) is equivalent
to the category consisting of pairs $(V,V')$, where $V\to X$ is a bundle
and where $V'$ is a subbundle of the restriction of $D$ to 
the $k$-th thickening $D_k$ of $D$ 
(if $\III\subset\OOO_X$ is the ideal sheaf defining $D$, then 
$D_k=\Spec(\OOO_X/\III^k)$). This can be proved by using the same
techniques as here. What seems more difficult is to understand the
category consisting of all the $\Gm$-equivariant bundles on $X_D$ in
terms of some category of bundles on $X$ with extra structure (defined
somehow along $D$). 
\end{remark}

The proof of the theorem is given in Subsections
\ref{firstthing} and \ref{secondthing}. 
The scheme of the proof is the following. First we prove
that for any pair of objects $W,W'\in\Veq(X_D)$ the map 
$$\mu_{W,W'}:\Mor_{\Veq(X_D)}(W,W')\to
\Mor_{\PPP(X,D)}(\mu(W),\mu(W'))$$
induced by the functor $\mu$ is a bijection.
In the second part of the
proof we will construct, for any $(V,V_1)\in\PPP(X,D)$, an
equivariant bundle $W\in\Veq(X_D)$ such that $\mu(W)$ is
isomorphic to $(V,V_1)$ (we use for that the result 
obtained in the first part).
That these two steps suffice to prove that $\mu$ induces an
equivalence of categories is assured by
Freyd's theorem (see Theorem 1.13 in Chapter 1 of \cite{GM}).

\begin{corollary}
Let $W\to X_D$ be a $\Gm$-equivariant bundle whose weights in $X_0$,
$X_1$ and $X_2$ are zero. Then, using the identification $X\simeq
X_2$ to write $p:X_D\to X_2$, we have a canonical isomorphism
$W\simeq p^*(W|_{X_2})$.
\label{peszerotrivial}
\end{corollary}
\begin{pf}
Let $W$ be such a vector bundle. Then $W\in\Veq(X_D)$, and $\mu(W)=(V,0)$ for
some vector bundle $V\to X$. On the other hand, if we define
$W_0:=p^*V$ and take on $W_0$ the trivial action of $\Gm$, it turns
out that $\mu(W_0)=(V,0)$ as well. Hence, by Theorem \ref{mainthm},
$W$ and $W_0$ are isomorphic as $\Gm$-vector bundles.
\end{pf}

\subsection{The map $\mu_{W,W'}$ is injective}
\label{firstthing}
Let $W,W'\in\Veq(X_D)$ be two $\Gm$-vector bundles.
Let $U=U^-(X_2)=\{z\in X_D\mid z_-\in X_2\}$.  
Since the only weights of $W$ at $X_2$ are
$0$, it turns out that if $z\in U$ then the map
$\rho^z_{W}:W_z\to W_{z_-}$ is an
isomorphism (see Lemma \ref{rhoiso}). Similarly, we have an isomophism
$\rho^z_{W'}:W'_z\to W_{z_-}$ whenever $z\in U$. Now, a morphism
$\phi\in\Mor_{\Veq(X_D)}(W,W')$ is by definition a $\Gm$-equivariant
map $\phi:W\to W'$, and equivariance means that for any $w\in W_z$ and
$\theta\in\Gm$ we have
$\phi_{\theta\cdot z}(\theta\cdot w)=\theta\cdot\phi_z(w),$
where $\phi_x:W_x\to W'_x$ is the restriction of $\phi$ to the fibres
over $x$.
Now suppose that $z\in U$ and make $\theta\to\infty$. We get at the
limit that $\phi_{z_-}\circ \rho^z_W=\rho^z_{W'}\circ\phi_z$, and since 
$\rho^z_{W'}$ is an isomorphism we can write
$$\phi_z=(\rho^z_{W'})^{-1}\circ\phi_{z_-}\circ \rho^z_W,$$
which means that $\phi|_U$ is determined by $\phi|_{X_2}$. And since
$U\subset X_D$ is dense, it follows that $\phi$ is also determined by
$\phi|_{X_2}$. Hence, $\mu_{W,W'}$ is injective.

\subsection{The map $\mu_{W,W'}$ is exhaustive}
\label{firstthing2}
Let $\phi:V=W|_{X_2}\to \tV=W'|_{X_2}$ be a morphism in $\PPP(X,D)$.
Let $U=U^-(X_2)$. 
To prove that $\mu_{W,W'}$ is exhaustive we use the following strategy:
first we define a $\Gm$-equivariant map $\psi:W|_{U}\to W'|_{U}$ which 
extends $\phi$, and then we prove that $\psi$ extends to a $\Gm$-equivariant
map $\psi:W\to W'$. It is in the second step that we use that $\phi$ is 
a morphism in $\PPP(X,D)$. 

We define $\psi$ as follows: for any $z\in U$,
\begin{equation}
\psi_z:=(\rho^z_{W'})^{-1}\circ\phi_{z_-}\circ\rho^z_W:W_z\to W'_z.
\end{equation}

Let $X_D^1=X_D\setminus X_1$. By Corollary \ref{corclosuresbdl}, in order
to prove that $\psi$ extends to the whole $X_D$ it
suffices to check that $\psi$ extends to $X_D^1$. 

Recall that $p:X_D\to X$ denotes the projection. 
We are going to prove that for any $z\in D\subset X$ there exists a
neighbourhood $z\in B\subset X$ such that $\psi|_{p^{-1}(B)\cap
U}$ extends to a map from $W|_{p^{-1}(B)\cap X_D^1}$ to
$W'|_{p^{-1}(B)\cap X_D^1}$  
This will be done by taking
metrics on $W|_{p^{-1}(B)}$ and $W'|_{p^{-1}(B)}$,
and by checking that the $L^{\infty}$
norm of $\psi|_{p^{-1}B\cap U}$ is bounded. The fact that $\psi$
extends will then follow from Riemann's extension theorem.
Note that since $p^{-1}(B)\cap U$ is dense in $p^{-1}(B)$ 
such extension must be unique, and this implies automatically
that the extensions of $\psi|_{p^{-1}B\cap U}$ for different choices
of $U$ patch together, and that the resulting extension of $\psi$ 
is $\Gm$-equivariant.

\subsubsection{}
\label{defnu}
Take a point $z\in D$ and coordinates $z_1,\dots,z_n$ in a neighbourhood
$B'\subset X$ of $z$ such that 
\begin{enumerate}
\item $z$ corresponds to $(0,\dots,0)$, 
\item $D$ is given by $z_1=0$, and 
\item $c=(z_1,\dots,z_n):B'\to\CC^n$ identifies
biholomorphically $B(0,1)\subset\CC^n$ with a neigbourhood $B\subset
B'$ of $z$ (here $B(a,r)$ is the ball of radius $r$
centered at $a\in\CC^n$). 
\end{enumerate}
Then there is a
$\Gm$-equivariant isomorphism $\nu:p^{-1}(B) \simeq
p_n^{-1}B(0,1)\subset H_n$ and we define $U':=U_n\cap p_n^{-1}B(0,1)$
(see Subsection \ref{defHn} for the definitions of $p_n:H_n\to\CC^n$
and $U_n\subset H_n$). 
Let us transport the vector bundles $W,W'$ to $p^{-1}B(0,1)$ by defining
$$Y=(\nu^{-1})^*W|_{p^{-1}B},\qquad Y'=(\nu^{-1})^*W'|_{p^{-1}B}.$$ 
By a slight abuse of notation,
$\psi: Y|_{U'} \to Y'|_{U'}$ will be the corresponding map. 
Take metrics on $Y$ and $Y'$ which are invariant under the
action of $S^1\subset\Gm$. 

\subsubsection{}
Let $\frac{1}{2}>\epsilon>\delta^2>0$ be small enough so that we have
$\Gm$-equivariant trivialisations as the ones given by Lemma
\ref{modellocalW}, and centered
around $z=(0,\dots,0,[0:0:1])$, 
of the restrictions of $Y,Y'$ to the open
subset $C=\xi^{-1}([-1-\delta,-1+\delta])\cap p_n^{-1}B(0,2\epsilon)$
of $H_n$. 
Let $S\subset C$ be the neighbourhood of $z$ defined as follows
$$S:=\{(t_1,\dots,t_n,[x:y:1])\mid \sum_{i=1}^n|t_i|^2\leq\epsilon^2,
\ 0\leq |x|,|y|\leq\delta,\ t_1=xy\}\subset p_n^{-1}B(0,1).$$

\begin{lemma}
The restriction of $\psi$ to $S\cap U'$ has bounded $L^{\infty}$
norm.
\label{boundS}
\end{lemma}
\begin{pf}
In all the proof $t$ will denote the $n$-uple $(t_1,\dots,t_n)$. 
Note that we have
$S\cap U'=\{(t,[x:y:1])\in S\mid x\neq 0\}.$
Let $Y_z=Y(0)\oplus Y(1)$ and $Y'_z=Y'(0)\oplus Y'(1)$ be the
decompositions in weights of the fibres of $Y$ and $Y'$ over $z$.
Then, since $S\subset C$, we have weakly $\Gm$-equivariant
trivialisations 
$$Y|_S\simeq S\times (Y(0)\oplus Y(1))\qquad\text{ and }\qquad
Y'|_S\simeq S\times (Y'(0)\oplus Y'(1)).$$
Using these trivialisations, the restriction of $\psi$ to $S\cap U'$
can be written in matrix form as
$$\psi=\left(\begin{array}{cc}
\psi_{11} & \psi_{01} \\ \psi_{10} & \psi_{00}\end{array}\right),$$
where we view $\psi_{ij}$ as a weakly $\Gm$-equivariant map from
$S\cap U'$ to $\Hom(Y(i),Y'(j))$. This means that if
$\theta\in\Gm$, $(t,[x:y:1])\in S$ and
$(t,[\theta x:\theta^{-1}y:1])\in S$, then
\begin{equation}
\psi_{ij}(t,[\theta x:\theta^{-1}y:1])
=\theta^{i-j}\psi_{ij}(t,[x:y:1]).
\label{equivariancia}
\end{equation}
Let $S_\delta=\{(t,[x:y:1])\in S\mid x=\delta\}.$

Let $S_{\delta,0}=\{(t,[x:y:1])\in S_\delta\mid y=0\}$, 
and observe that $S_{\delta,0}\subset U^+(z)$.
Since the map $\phi$ which was used
to construct $\psi$ belonged to $\PPP(X,D)$, we deduce that 
$\psi(W^{+,1}(X_0))\subset {W'}^{+,1}(X_0)$.
The restrictions to $S_{\delta,0}$ 
of the pull backs $\nu^*W^{+,1}(X_0)$ and $\nu^*{W'}^{+,1}(X_0)$ correspond, 
using the trivialisations of $Y$ and $Y'$, 
to $S_{\delta,0}\times Y(1)$ and $S_{\delta,0}\times Y'(1)$
respectively, and from this it follows that
\begin{equation}
\psi_{10}|_{S_{\delta,0}}=0.
\label{sanulla}
\end{equation}
It now follows from (\ref{sanulla}) and the facts that
$\psi$ is smooth and $S_\delta\subset S\cap U'$ is compact, 
that there is a constant $K$ such that for any 
$(t,[\delta:y:1])\in S_\delta$ we have 
\begin{align}
|\psi_{10}(t,[\delta:y:1])| & \leq K|y| \label{fita1} \\
|\psi_{ij}(t,[\delta:y:1])| & \leq K,   \label{fita2}
\end{align}
for any $(i,j)\subset\{0,1\}^2$. 
Now, combining (\ref{equivariancia}) with (\ref{fita1}) and
(\ref{fita2}) we deduce that $\psi|_{S\cap U'}$ has bounded
$L^{\infty}$ norm. 
Indeed, let $(t,[x:y:1])\in S\cap U'$. Putting 
$\theta=\delta x^{-1}$ in (\ref{equivariancia}) and taking norm we
get, using (\ref{fita1}):
$$|\psi_{10}(t,[x:y:1])| =
|\delta x^{-1} \psi_{10}(t,[\delta:yx\delta^{-1}:1])| 
\leq |\delta x^{-1} K yx\delta^{-1}| 
= |Ky|\leq K\delta.$$
(This makes sense because 
$(t,[\delta:yx\delta^{-1}:1])\in S_\delta$, since 
$|yx\delta^{-1}|\leq\delta$.) Finally, if $(1,0)\neq (i,j)\in\{0,1\}^2$ and 
$(t,[x:y:1])\in S\cap U'$ we obtain, using
(\ref{equivariancia}) with $\theta=\delta x^{-1}$ and (\ref{fita2}):
$$|\psi_{ij}(t,[x:y:1])| =
|(\delta x^{-1})^{i-j}\psi_{ij}(t,[\delta:yx\delta^{-1}:1])| 
\leq|\psi_{ij}(t,[\delta:yx\delta^{-1}:1])| 
\leq K,$$
since then $i-j\leq 0$, so $|x|\leq\delta$ implies $|\delta
x^{-1}|^{i-j}\leq 1$. 
\end{pf}

\subsubsection{}
\label{lastthing}
Let $0<\delta'<\delta$ and $0<\epsilon'<\epsilon$ be small enough so
that $C'=\xi^{-1}([-1-\delta',-1+\delta'])\cap
p_n^{-1}B(0,\epsilon')\subset S$.

\begin{lemma}
For any $0>\alpha>-1+\delta'$ 
the restriction of $\psi$ to $\xi^{-1}([-2,\alpha])\cap
p_n^{-1}B(0,\epsilon')\cap U$ has bounded 
$L^{\infty}$ norm.
\end{lemma}
\begin{pf}
Let $C_1:=\xi^{-1}([-2,-1-\delta'])\cap
p_n^{-1}B(0,\epsilon')$.
Since $\ov{C_1}\subset U$, 
it follows from compactness that
$|(\psi|_{C_1})|_{L^{\infty}}<\infty$. 
By assumption,
$C_2:=C'\cap U\subset S$, so Lemma \ref{boundS} implies that
$|(\psi|_{C_2})|_{L^{\infty}}<\infty$. It remains to prove that
$|(\psi|_{C_3})|_{L^{\infty}}<\infty$, where 
$C_3=\xi^{-1}([-1+\delta',\alpha])\cap p_n^{-1}B(0,\epsilon')\cap U$.
Given $z\in\ov{C_3}$, let
$\lambda(z)\in\RR$ be defined by the condition 
$\xi(\lambda(z)\cdot z)=-1+\delta$ (this is well defined because the
elements of $\ov{C_3}$ have trivial stabiliser). By compactness of
$\ov{C_3}$ it follows
that $\lambda:C_3\to\RR$ is bounded. Finally, since
$\xi^{-1}(-1+\delta)\cap C_3\subset C_2$ and we already know that 
$|(\psi|_{C_2})|_{L^{\infty}}<\infty$, it follows from
the $\Gm$-equivariance of $\psi$ and the boundedness of $\lambda$
that $|(\psi|_{C_3})|_{L^{\infty}}<\infty$. The lemma follows now because 
$$\xi^{-1}([-2,\alpha])\cap
p_n^{-1}B(0,\epsilon')\cap U= C_1\cup C_2\cup C_3.$$
\end{pf}

This lemma implies, by making $\alpha\to 0$, 
that the section $\psi$ extends to $X_D^1$, and this finishes the
proof of the exhaustivity of $\mu_{W,W'}$.

\subsection{Constructing objects of $\Veq(X,D)$ from objects of $\PPP(X,D)$}
\label{secondthing}

Let $\LLL\to X_D$ be the line bundle associated to the divisor 
$\Delta^-$, with the lift of the
action of $\Gm$ whose existence is granted by Lemma \ref{useful2}.
By construction it follows that there is a $\Gm$-equivariant
isomorphism $\LLL|_{\Delta^-}\simeq N^-$, so from Lemma \ref{weightN2}
we deduce that $\LLL$ is an object of $\Veq(X_D)$ (that the weight of
the restriction $\LLL|_{X_2}$ is zero follows from $\Delta^-\cap
X_2=\emptyset$, see Lemma \ref{useful2}).
Finally, it follows from the construction that there is a nowhere zero
section $\psi\in H^0(X_2;\LLL|_{X_2})$ (again, because $X_2\cap\Delta^-=
\emptyset$). 

Let us take an object $(V,V_1)\in\PPP(X,D)$. Our aim is to find an
equivariant bundle $W\in\Veq(X_D)$ such that $\mu(W)\simeq (V,V_1)$.
Let $E=E_0\oplus E_1$ be a direct sum of vector spaces satisfying $\dim
E=\rk V$ and $\dim E_1=\rk V_1$.
Let $\{A_q\}_{q\in I}$ be a family of open subsets of $X$, which cover
$X$, and such that for any $q\in I$
there is a vector bundle isomorphism $\phi_q:V|_{A_q}\to
A_q\times E$ such that $\phi_q(V_1|_{A_q\cap D})=(A_q\cap D)\times E_1$.
Let us define for any $q,q'$ the following sets:
$$D_q=A_q\cap D,\qquad A_{q,q'}=A_q\cap A_{q'},\qquad 
D_{q,q'}=A_q\cap A_{q'}\cap D.$$ 
The isomorphism class of
$(V,V_1)$ is completely determined by the transition functions
$$\{\phi_{q,q'}=\phi_{q'}|_{A_q\cap A_{q'}}
\circ(\phi_q|_{A_q\cap A_{q'}})^{-1} \mid q,q'\in I\}.$$

Note that the covering $\{A_q\}$ of $X$ induces a covering
$\{(A_q)_{D_q}\}$ of $X_D$. To save on typing, we will denote by
$p$ the restriction of the projection map $p:X_D\to X$ on any
$(A_q)_{D_q}$.
We define for any $q\in I$ the following $\Gm$-equivariant bundle on
$(A_q)_{D_q}$: 
$$M_q:=\underline{E_0}\oplus \underline{E_1}\otimes
\LLL|_{(A_q)_{D_q}}$$ 
(here $\underline{E_0}$ and $\underline{E_1}$
denote the trivial bundles on $(A_q)_{D_q}$ 
with fibres $E_0$, $E_1$ respectively, and with the 
trivial lift of the action of $\Gm$).
Then $M_q$ is an object of $\Veq((A_q)_{D_q})$.
To get an equivariant bundle on $X_D$ by patching the bundles $M_q$
we have to specify transition functions 
$$\{\Phi_{q,q'}:M_q|_{(A_q)_{D_q}\cap (A_{q'})_{D_{q'}}}
\to M_{q'}|_{(A_q)_{D_q}\cap (A_{q'})_{D_{q'}}} 
\mid q,q'\in I\}.$$
The section $\psi$ defined above provides, for any $q$, and isomorphism 
$$\eta'_q:A_q\times E\to M_q|_{X_2\cap (A_q)_{D_q}}$$
(here we use the identification $X_2\simeq X$ to identify
$A_q\simeq X_2\cap (A_q)_{D_q}$), which induces an isomorphism
of objects in $\PPP(A_q,D_q)$
$$\eta_q:(A_q\times E,D_q\times E_1)\to\mu(M_q).$$
Now, to define $\Phi_{q,q'}$ we will use the fact that the map 
$$\mu_{q,q'}:\Mor_{\Veq((A_q)_{D_q}\cap (A_{q'})_{D_{q'}})}
(M_q,M_{q'})\to
\Mor_{\PPP(A_{q,q'},D_{q,q'})}(\mu(M_q,M_{q'})$$
induced by the functor $\mu$
is an isomorphism (this has been proved in Sections \ref{firstthing}
and \ref{firstthing2}). Namely, we set:
$$\Phi_{q,q'}:=\mu_{q,q'}^{-1}
(\eta_{q'}\circ\phi_{q,q'}\circ\eta_{q}^{-1}).$$
Then the functions $\{\Phi_{q,q'}\}$ are all $\Gm$-equivariant and
they satisfy the cocycle condition (because the functions
$\{\phi_{q,q'}\}$ satisfy it), so they define an equivariant bundle
$W\to X_D$. The isomorphisms $\{\eta_q\}$ and $\{\phi_q\}$ give
an isomorphism between $\mu(W)$ and $(V,V_1)$, and we are done.

\subsection{An equivariant version of Theorem \ref{mainthm}}
\label{equivver}
Suppose that a group $G$ acts on $X$ preserving the divisor
$D$. By Lemma \ref{useful} such an action lifts to 
$X_D$ and the resulting action commutes with the action of
$\Gm$. Furthermore, $X_0\cup X_1\cup X_2\subset X_D^{G}$.
Consider the following categories (which are the $G$-equivariant
versions of the categories $\Veq(X_D)$ and $\PPP(X,D)$):
\begin{enumerate}
\item $\Veq^{G}(X_D)$ is the category whose objects are
  $G\times\Gm$-equivariant vector bundles $W\to X_D$ which,
  considered as a $\Gm$-equivariant bundle, belong to $\Veq(X_D)$, and
  whose morphisms are the $G\times\Gm$-equivariant morphisms of
  vector bundles.
\item $\PPP^{G}(X,D)$ is the category whose objects are pairs
  $(V,V_1)$, where $V\to X$ is a $G$-equivariant vector bundle and
  $V_1\to D$ is a $G$-invariant subbundle of $V|_D$, and whose
  morphisms are the $G$-equivariant morphisms of bundles which
  preserve the subbundles.
\end{enumerate}
Note that we have obvious functors
$f_{\VVV}:\Veq^{G}(X_D)\to\Veq(X_D)$ and
$f_{\PPP}:\PPP^{G}(X,D)\to\PPP(X,D)$ which forget the action of
$G$. We now construct a functor
$\mu^{G}:\Veq^{G}(X_D)\to\PPP^{G}(X,D)$ as follows. The
action of $\mu^{G}$ on objects maps and equivariant bundle
$W\to\Veq^{G}(X_D)$ to the pair
$\mu^{G}(W)=(V,V_1)\in\PPP^{G}(X,D)$ satisfying: (1)
$f_{\PPP}(V,V_1)=\mu(f_{\VVV}W)$ and (2) the lift of $G$ to $V$ is
the one obtained by identifying $V=W|_{X_2}$ (recall that we have
canonically $X_2\simeq X$). Finally, the action of $\mu^{G}$ on
morphisms is by taking the restriction to $X_2$ (exactly like that of
$\mu$). 

\begin{theorem}
The functor $\mu^{G}$ induces an equivalence of categories.
\label{eqmainthm}
\end{theorem}
\begin{pf}
This follows essentially from Theorem \ref{mainthm} and the
commutativity of diagram (\ref{commutamu}) applied to the
biholomorphisms $f\in\Bihol(X)$ given by the action of the elements of
$G$ on $X$. 
The only thing which might not be clear is the fact that any object in
$\PPP^{G}(X,D)$ is isomorphic to the image by $\mu^{G}$ of
some object in $\Veq^{G}(X_D)$. Let us clarify this point.
Giving a lift of the $G$ action on $X$ to a vector bundle $V\to X$ is
the same thing as giving a set of isomorphisms
$$I=\{i_g:V\to l_g^*V\mid g\in G\}$$
(where $l_g:X\to X$ denotes the biholomorphism induced by the action of
$g\in G$) which satisfy certain cocycle condition (let us call it
$G$-condition). 
If $(V,V_1)\in\PPP^{G}(X,D)$, the isomorphisms in $I$ are also
compatible with $V_1$. Hence, if we have $W\in\Veq(X_D)$ and an
isomorphism $\mu(W)\simeq V$, we can use Theorem \ref{mainthm} to
lift the family of isomorphisms $I$ to a family of
$\Gm$-equivariant isomorphisms
$$I'=\{i_g:W\to {l'_g}^*W\mid g\in G\},$$
denoting $l'_g:X_D\to X_D$ the action of $g\in G$. 
Finally, again by Theorem \ref{mainthm}, the elements in $I'$ satisfy
the $G$-condition, precisely because the elements in $I$ satisfy it.
\end{pf}

In the following lemmae, recall that $N^-\to\Delta^-$ and
$N^+\to\Delta^+$ are the normal bundles of the inclusions
$\Delta^-\subset X_D$ and $\Delta^+\subset X_D$.
By an abuse of notation, we will denote by $L\to D$ be the normal
bundle of the inclusion $D\subset X$. 

\begin{lemma}
Suppose that $G$ is reductive.
Let $a\in X^{G}$, $W\in\Veq^{G}(X_D)$ and
$(V,V_1):=\mu^{G}(W)$. Let $b\in p^{-1}(a)\cap X_D^{G}$. If
$b\notin\Delta^-$ (resp. if $b\in\Delta^-$) then there is a
$G$-equivariant isomorphism of fibres $W_b\simeq V_a$
(resp. $W_b\simeq (V_1)_a\otimes N^-_b\oplus (V_a/(V_1)_a)$).
\label{cincdos}
\end{lemma}
\begin{pf}
Since by assumption $G$ is reductive and everything is
holomorphic, we can use Weyl's unitary trick and work with a maximal
compact subgroup of $G$, which by an abuse of notation we will
denote, only in this proof, by the same symbol $G$. 
If $a\notin D$ then the result follows easily from Lemma \ref{rhoiso}.
Now assume that $a\in D$. Using the same idea as in the proof of Lemma
\ref{modellocalW} and up to shrinking $X$ to a small
$G$-invariant neighbourhood of $a$ we may assume that there is a
$G$-equivariant trivialisation 
$$V\simeq (X\times (V_1)_a)\oplus(X\times V_a/(V_1)_a)$$ inducing a
trivialisation $V_1\simeq D\times (V_1)_a$. 
Let $\LLL\to X_D$ be the line bundle corresponding to the divisor
$\Delta^-$ together with the lift of the $G\times\Gm$ action
induced by Lemma \ref{useful2}.
Then $$W_0:=(X_D\times V_a')\otimes\LLL\oplus(X_D\times(V_a/(V_1)_a))$$
with the corresponding lifts of the $G\times\Gm$ action satisfies
$\mu^{G}(W_0)=(V,V_1)$, so by Theorem \ref{eqmainthm} we have a
$G\times\Gm$-equivariant isomorphism $W_0\simeq W$. Now the
result follows by using the second part of Lemma \ref{useful2} to
identify $\LLL_b\simeq N^-_b$ as representations of $G$.
\end{pf}

\begin{lemma}
For any $b\in X_D^{G}\cap\Delta^-$ there is a $G$-equivariant
isomorphism $N^-_b\simeq L_{p(b)}$.
\label{cinctres}
\end{lemma}
\begin{pf}
Let $a=p(b)$ and identify
$\PP(\CC^3)\simeq\PP(L\oplus\uC\oplus\uC)_a$. Since
$p^{-1}(a)\cap\Delta^-\cap X_D^{G}$ is either
$\{[0:1:0]\}\cup\{[0:0:1]\}$ or $\{[0:x:y]\mid
(x,y)\in\CC^2\setminus(0,0)\}$ (depending on whether the representation 
$L_a$ of $G$ is trivial or not), it suffices to prove the
result for $b=[0:1:0]\in X_1$ or $b=[0:0:1]\in X_2$. 
For the first case, observe that the restriction of $N^-$ to
$X_1\cap\Delta^-$ is isomorphic to $L$ (because $X_1$ and $\Delta^-$
intersect transversely along a submanifold which can be canonically
identified --- using $p$ --- with $D$). If $b=[0:0:1]$, then observe
that the map $f:L_a\to X_D$ which sends $x\in L_a$ to $[x:0:1]$ is
equivariant and transverse to $\Delta^-$ at $b$. Hence it gives an
equivariant identification of $L_a$ with $N^-_b$. On the other hand, by
Lemma \ref{useful2} the representation $L_a$ of $G$ is isomorphic
to $L_a$, and we are done.
\end{pf}

\begin{lemma}
Let $b\in X_D^{G}\cap \Delta^+$. If $b\in X_0$ then the
representation $N^+_b$ of $G$ is the trivial one, and if $b\in X_2$
then the representation $N^+_b$ is isomorphic to $L_{b(p)}$.
\end{lemma}
\begin{pf}
Exactly the same as that of the previous lemma.
\end{pf}

\section{Parabolic structures over a smooth divisor}
\label{divllis}
In this section $X$ will be a manifold and $D\subset X$ a smooth
divisor. Let us fix an integer $r\geq 1$.
Let $\PPP(X,D,r)$ be the category defined as follows:
\begin{enumerate}
\item The objects of $\PPP(X,D,r)$ are pairs
  $(V,\VVV)$ consisting of a vector bundle $V\to X$ and a
  filtration $\VVV$ of vector bundles over $D$: 
  $$\VVV=(V_1\subset\dots\subset V_r\subset V|_D)$$ 
  (note that the inclusions need not be strict).
\item The morphisms in $\PPP(X,D,r)$ between two objects 
  $(V,\VVV)$ and $(V',\VVV')$ are the morphisms of
  vector bundles $\phi:V\to V'$ whose restriction to $D$ respects
  the filtrations $\VVV$ and $\VVV'$, i.e., for
  any $1\leq j\leq r$, $\phi|_D(V_j)\subset V'_j$.
\end{enumerate}
Our aim in this section is to obtain, by applying recursively the
construction of Section \ref{construccioXD}, a $(\Gm)^r$-manifold
$X(D,r)$ which fibers 
over $X$, and to prove that the category $\PPP(X,D,r)$ is equivalent to
a full subcategory of the category of $(\Gm)^r$-equivariant bundles
over $X(D,r)$. 

\subsection{Notations and definitions}

\subsubsection{Weights of $(\Gm)^r$}
\label{defweights}
The objects of the subcategory of equivariant vector bundles which
will ultimatelly be equivalent to $\PPP(X,D,r)$ will be those which
satisfy a certain restriction on the weights of a action of
$(\Gm)^r)$ (to be defined below). In order to be able to specify this
restriction we need to introduce some notation on weights.

We identify the characters of $(\Gm)^r$ with $\ZZ^r$ by assigning to
$a=(a_1,\dots,a_r)\in\ZZ^r$ the character
$$\chi_a:(\Gm)^r\ni(\theta_1,\dots,\theta_r)\mapsto
\theta_1^{a_1}\cdot\dots\cdot\theta_r^{a_r}\in\Gm.$$
If $r<s$ we map $\ZZ^r\to\ZZ^s$ by sending $(a_1,\dots,a_r)$ to
$(a_1,\dots,a_r,0,\dots,0)$, and in this way we will sometimes
implicitely view the set of characters of $(\Gm)^r$ as subset of 
the characters of
$(\Gm)^s$. In particular, $0\in\ZZ^r$ will denote $(0,\dots,0)$.

For any $1\leq j\leq r$ let $\pi_j:\ZZ^r\to\ZZ$ be the projection to the 
$j$-th factor. Let $e_j$ be the character corresponding
to the projection $(\Gm)^r\to\Gm$ to the $j$-th factor. Note that 
$e_1,\dots,e_r$ is the canonical basis of $\ZZ^r$ and that
$\pi_i(e_j)=\delta_{ij}$. We also define for
any $1\leq j\leq r$
$$f_j:=e_j+2e_{j-1}+\sum_{i\geq 2}e_{j-i},$$ 
(where we understand that $e_k=0$ whenever $k\leq 0$).

Let $\Pi_{j+1}:\ZZ^{j+1}\to\ZZ^j$ be defined as
$$\Pi_{j+1}(f):=\left\{\begin{array}{ll}
f & \text{if }\pi_{j+1}(f)=0 \\
f-\pi_{j+1}(f)-e_{j+1}-e_j+e_{j+1} & \text{if }\pi_{j+1}(f)\neq 0.
\end{array}\right.$$
The following equalities follow immediately from the definition:
\begin{equation}
\Pi_{j+1}(f_{j+1})=f_j\qquad\text{ and }\qquad
\Pi_{j+1}(f_i)=f_i\text{ for any $i\leq j$}.
\label{Quatre}
\end{equation}

\subsubsection{The iterated construction}
\label{iterconstr}
We define recursively a sequence of manifolds $$Y_0,Y_1,\dots,Y_r$$
with projections $p_j:Y_j\to Y_{j-1}$ and a smooth divisor
$\Delta^+_j\subset Y_j$ for any $1\leq j\leq r$ as follows.
We first set $Y_0:=X$ and $\Delta^+_0:=D$. If $0\leq j<r$, we
apply the construction of Subsection \ref{construccioXD} to define
$$p_{j+1}:Y_{j+1}:=(Y_j)_{\Delta^+_j}\to Y_j.$$
By construction the manifold $Y_{j+1}$ carries an action of $G_j:=\Gm$,
whose fixed point locus is the disjoint union of two copies of $Y_j$,
denoted by $Y_{j,1},Y_{j,2}$, and a copy of $\Delta^+_j$,
denoted by $Y_{j,0}$. With this in mind, we define
$$\Delta^+_{j+1}:=\overline{U^+(Y_{j,0})}\qquad
\text{ and }\qquad
\Delta^-_{j+1}:=\overline{U^-(Y_{j,0})}.$$
Both $\Delta^+_{j+1}$ and $\Delta^-_{j+1}$ are smooth divisors in of $Y_{j+1}$.

Let $i\in\{0,1,2\}$. We define $Y_1[i]:=Y_{0,i}$ and, for any $1\leq j<r$, 
we set $Y_{j+1}[i]:=Y_{j,i}\cap p_{j+1}^{-1}(Y_j[i])\subset Y_{j+1}$.

If $1\leq j<r$, we may use Lemma \ref{useful} to lift the action of
$G_j$ on $Y_j$ to an action on $Y_{j+1}$, which commutes with the
action of $G_{j+1}$. And, using recursion, we get, for any $1\leq
j\leq r$, a natural action of $G(j)=G_1\times\dots\times G_j$ on
$Y_j$. It is easy to check that for any $i\in\{0,1,2\}$ the
submanifold $Y_j[i]$ belongs to the fixed point set of
the action of $G(j)$.

\begin{definition}
We define $X(D,r)$ to be $Y_r$ and, for any $i\in\{0,1,2\}$, 
$X(D,r)[i]:=Y_r[i]$.
\end{definition}

By construction we have a tower of manifolds 
\begin{equation}
X(D,r)=Y_r\stackrel{p_r}{\longrightarrow}
Y_{r-1}\stackrel{p_{r-1}}{\longrightarrow}\dots
\stackrel{p_1}{\longrightarrow}
Y_0=X
\label{tow1}
\end{equation}
and, identifying for any $0\leq j<r$ the submanifold
$Y_j$ with $Y_{j,2}\subset Y_{j+1}$ we get a chain of inclusions
\begin{equation}
\label{inclXj}
X=Y_0\subset Y_1\subset\dots\subset Y_{r-1}\subset Y_r=X(D,r).
\end{equation}
These inclusions induce a chain of equalities
\begin{equation}
X=Y_1[2]=Y_2[2]=\dots=Y_r[2].
\label{duesestr}
\end{equation}

\begin{lemma}
For any $i\in\{0,1,2\}$ the following is a strong inclusion
(see Subsection \ref{def:stronginclusion}):
$Y_{j+1}[i]\subset Y_{j,i}.$ 
\label{strincl}
\end{lemma}
\begin{pf}
Easy from the definitions.
\end{pf}

\subsubsection{Categories of equivariant bundles}

Let $\Veqr(X(D,r))$ be the category defined as follows.
\begin{enumerate}
\item The objects of $\Veqr(X(D,r))$ are $G(r)$-equivariant bundles
  $W\to X(D,r)$ whose only $G(r)$-weight on $X(D,r)[1]$ and $X(D,r)[2]$ is 
  zero and whose $G(r)$-weights on $X(D,r)[0]$ belong to
  $\{f_1,\dots,f_r\}$.
\item The morphisms in $\Veqr(X(D,r))$ between to objects $W$ and $W'$
  are the $G(r)$-equivariant maps of vector bundles $\psi:W\to W'$.
\end{enumerate}
Our aim is to construct a functor $\mu(r):\Veqr(X(D,r))\to\PPP(X,D,r)$
inducing an equivalence of categories.
For that we will construct a set of auxiliar categories
$\VP_0,\dots,\VP_r$ making up a bridge from $\Veqr(X(D,r))$ 
to $\PPP(X,D,r)$.

Let $0\leq j\leq r$. Define the category $\VP_{j}$ as follows.
\begin{enumerate}
\item The objects of $\VP_{j}$ are pairs $(W^j,\WWW^j)$
consisting of a $G(j)$-equivariant bundle
$W^j\to Y_j$ and a filtration of bundles over $\Delta^+_j$:
$$\WWW^j=(W^j_1\subset\dots\subset W^j_{r-j}\subset W^j|_{\Delta^+_j}),$$
subject to the following constraints: the only $G(j)$-weight of
$W^j$ on $Y_j[1]$ and $Y_j[2]$ is zero, the $G(j)$-weights
of $W^j$ on $Y_j[0]$ belong to $\{f_1,\dots,f_j\}$ and, for any
$1\leq s\leq r-j$, the only $G(j)$-weight of $W^j_s$ on $Y_j[0]$ is $f_j$
(in fact this condition holds for any $1\leq s\leq r-j$ if and only
if it holds for $s=r-j$).

\item The morphisms in $\VP_{j}$ between two objects
$(W^j,\WWW^j)$ and $({W'}^j,{\WWW'}^j)$ are the $G(j)$-equivariant maps of
vector bundles $\phi:W^j\to {W'}^j$ whose restriction to $\Delta^+_j$
respects the filtrations $\WWW^j$ and ${W'}^j$.
\end{enumerate}
Observe that $\VP_{r}=\VVV_{G(r)}(X(D,r))$ and
$\VP_{0}=\PPP(X,D,r)$. In the next subsection we will 
construct for any $0\leq j<r$
a functor $\mu_{j+1}:\VP_{j+1}\to\VP_{j}$, and we will define
$\mu(r)$ to be the composition
$\mu(r):=\mu_1\circ\mu_1\circ\dots\circ\mu_r$. Finally, we will
show that each $\mu_j$ induces an equivalence of categories.
Hence, $\mu(r)$ gives also an equivalence of categories.

\subsection{The functors $\mu_{j+1}$}
\label{defmuj}

Let us fix some $0\leq j<r$. Our aim here is to define a functor
$$\mu_{j+1}:\VP_{j+1}\to\VP_{j}.$$

\subsubsection{}
Fist of all, we define the action of $\mu_{j+1}$ on morphisms 
to be the restriction to $Y_j\subset Y_{j+1}$.

We now define $\mu_{j+1}$ acting on objects. Let 
$(W^{j+1},\WWW^{j+1})\in\VP_{j+1}$, where
$\WWW^{j+1}=(W^{j+1}_1\subset\dots\subset W^{j+1}_{r-j-1})$. 
By definition, $W^{j+1}$ is a vector bundle over 
$Y_{j+1}=(Y_j)_{\Delta^+_j}$. 
We can now use the functor $\mu^{G(j)}$ defined in Subsection
\ref{equivver} to define
$$(W^j,W^j_{r-j}):=\mu^{G(j)}(W^{j+1}).$$
On the other hand, recall that we have an inclusion $Y_j\subset
Y_{j+1}$ by identifying $Y_j=Y_{j,2}\subset Y_{j+1}$. 
In this way we get an
identification between $\Delta^+_j$ and $\Delta^+_{j+1}\cap Y_{j,2}$.
With this in mind, we define, for any $1\leq s\leq r-j-1$, 
$$W^j_s:=W^{j+1}_s|_{\Delta^+_{j+1}\cap Y_{j,2}},$$
and we set $\WWW^j:=(0\subset W^j_1\subset\dots\subset W^j_{r-j})$.
Then, we define
$$\mu_{j+1}(W^{j+1},\WWW^{j+1}):=(W^j,\WWW^j).$$

\subsubsection{}
Let us check that $(W^j,\WWW^j)$ is indeed an object of $\VP_{j}$.
To begin with, observe that, since both $Y_{j,2}$ and
$\Delta^+_{j+1}$ are $G(j)$ invariant submanifolds of $Y_{j+1}$, all the
bundles $W^j,W^j_1,\dots,W^j_{r-j}$ inherit an action of $G(j)$. 

We also need to check that $W^j_{r-j-1}\subset W^j_{r-j}$.
By assumption, for any $1\leq s\leq r-j-1$, 
the only $G(j+1)$-weight of $W^j_s$ on $Y_{j+1}[0]$ is $f_{j+1}$.
This implies that the only $G_{j+1}$-weight of $W^j_s$ on 
$Y_{j+1}[0]$ is $1$.  Indeed, since
$Y_{j+1}[0]\subset Y_{j,0}$ is a strong inclusion (see Lemma
\ref{strincl}) and $Y_{j,0}\subset Y_{j+1}^{G_{j+1}}$, we know that 
the only $G_{j+1}$-weight of $W^j_s$ on $Y_{j,0}$ is $1$.
Now, Lemma \ref{inclsubbdl} implies that $W^j_{r-j-1}\subset W^j_{r-j}$.

In remains to check that the weights of $W^j$ restricted to $Y_j[i]$
for $i=0,1,2$ are the right ones. This is proved in the following lemma.

\begin{lemma}
Let $(W^j,W^j_{r-j})=\mu^{G(j)}(W^{j+1}).$ Then the weights of 
$G(j)$ acting on the restriction of $W^j$ to 
$Y_j[1]$ (resp. $Y_j[2]$, $Y_j[0]$) are $0$ (resp. $0$, contained
in $\{f_1,\dots,f_j\}$). 
Furthermore, $\chi_{G(j)}(W^j_{r-j}|_{Y_j[0]})=f_j$.
\end{lemma}
\begin{pf}
(a) Observe first that $G(j)$ acts trivially on 
$p^{-1}_{j+1}Y_j[1]\subset Y_{j+1}$. Indeed, the action of $G(j)$ on
$Y_j[1]\subset Y_j$ is trivial. And, since $Y_{j,1}\cap\Delta^+_j=\emptyset$,
the lift of the action of $G(j)$ to the line bundle $L_j\to Y_j$ 
corresponding to the divisor $\Delta^*(j)\subset Y_j$ is trivial
(see Lemma \ref{useful2}).
On the other hand, $Y_{j+1}[1]\subset p^{-1}_{j+1}Y_j[1]$ is a strong 
inclusion. Hence, since by hypothesis
$\chi_{G(j)}(W^{j+1}|_{Y_{j+1}[1]})=\{0\}$, it follows that
$\chi_{G(j)}(W^{j+1}|_{p^{-1}_{j+1}Y_j[1]})=\{0\}.$
Finally, by definition
$W^{j}|_{Y_j[1]}=W^{j+1}|_{p^{-1}_{j+1}Y_j[1]\cap Y_{j,2}}$
as $G(j)$-equivariant bundles, so that 
$\chi_{G(j)}(W^j|_{Y_j[1]})=0$.

(b) By (\ref{duesestr}) and the definition of $W^j$ we have
$\chi_{G(j)}(W^j|_{Y_j[2]})=\chi_{G(j)}(W^{j+1}|_{Y_{j+1}[2]})=0,$
where the second equality follows from our hypothesis.

(c) For any $1\leq s\leq r$, let $N^-_j\to\Delta^-_j$ and
$N^+_j\to\Delta^+_j$ be the normal bundles of the inclusions 
$\Delta^-_j\subset Y_j$ and $\Delta^+_j\subset Y_j$ respectively.
Note that $N^+_j=L_j|_{\Delta^+_j}$, where $L_j\to Y_j$ is
the line bundle corresponding to $\Delta^+_j\subset Y_j$.
We will now compute the weight of each group $G_i$
acting on $N^-_j|_{Y_j[0]}$.
\begin{enumerate}
\item Suppose that $i<j-1$. Then
$\chi_{G_i}(N^-_j|_{Y_j[0]}) = \chi_{G_i}(N^+_{j-1}|_{Y_{j-1}[0]})=0$,
where the first equality follows from by Lemma \ref{cinctres},
and the second one from $Y_{j-1}[0]\cap Y_{j-2,2}=\emptyset$.
\item The case $i=j-1$. We have
$\chi_{G_{j-1}}(N^-_j|_{Y_j[0]}) = 
\chi_{G_{j-1}}(N^+_{j-1}|_{Y_{j-1}[0]})= -1$,
where the first equality follows from Lemma \ref{cinctres},
and the second one from
$\chi_{G_{j-1}}(N^+_{j-1}|_{Y_{j-2,0}})=\{-1\}$ (Lemma \ref{weightN})
and $Y_{j-1}[0]\subset Y_{j-2,0}$.
\item The case $i=j$. We have
$\chi_{G_{j}}(N^-_j|_{Y_j[0]}) = 1$
since, by Lemma \ref{weightN2},
$\chi_{G_{j}}(N^-_j|_{Y_{j-1,0}})=\{1\},$ and $Y_j[0]\subset Y_{j-1,0}$.
\end{enumerate}
We have thus obtained:
\begin{equation}
\chi_{G(j)}(N^-_j|_{Y_j[0]})=-e_{j-1}+e_j.
\label{Cinc}
\end{equation}
Let now $a\in Y_j[0]$ and $b:=p^{-1}_{j+1}(a)\cap Y_{j,0}\in
Y_{j+1}[0]$. Let $\chi_b\subset\ZZ^{j+1}$ 
(resp. $\chi_a'\subset\chi_a\subset\ZZ^j$) 
be the weights appearing in the decomposition of $(W^{j+1})_b$
(resp. $(W^j_{r-j})_a\subset (W^j)_a$) in irreps of
$G(j+1)$ (resp. $G(j)$).
Combining Lemma \ref{cincdos} with (\ref{Cinc}) we deduce that
$$\chi_a=\Phi_{j+1}(\chi_b)\qquad
\text{ and }
\qquad \chi_a'=\{\Pi_{j+1}(\chi)\mid
\chi\in\chi_b,\ \pi_{j+1}(\chi)=1\}.$$
Now, by assumption
$\chi_b\subset\{f_1,\dots,f_{j+1}\},$
so by (\ref{Quatre}) we deduce that $\chi_a\subset\{f_1,\dots,f_{j}\}$.
Finally, if $f\in\{f_1,\dots,f_{j+1}\}$ then 
$\pi_{j+1}(f)=1$ if and only if $f=f_{j+1}$. Consequently, again
by (\ref{Quatre}), $\chi_a'=\{f_j\}$.
This finishes the proof.
\end{pf}

\begin{lemma}
The functor $\mu_{j+1}$ induces an equivalence of categories.
\label{mainthmj}
\end{lemma}
\begin{pf}
The proof is completely analogous to that of Theorems
\ref{mainthm} and \ref{eqmainthm}. The key step is to prove that
$\mu_j$ induces bijections when acting on morphisms, and this
can be done following exactly the same steps as in Subsections
\ref{firstthing} and \ref{firstthing2}.
\end{pf}

\subsubsection{} We sum up in the following theorem all the results which
we have obtained so far.

\begin{theorem} 
\label{mainthmr}
For any complex manifold $X$, any smooth divisor
$D\subset X$, and integer $r\geq 1$, there exists 
\begin{enumerate}
\item a manifold $X(D,r)$ acted on by $G(r)=(\Gm)^r$, 
\item an invariant projection $\pi:X(D,r)\to X$ with a section
$\sigma:X\to X(D,r)$, 
\item a full subcategory $\Veqr(X(D,r))$ of the category
of $G(r)$-equivariant vector bundles on $X(D,r)$, and 
\item a functor $\mu(r):\Veqr(X(D,r))\to\PPP(X,D,r)$ inducing
an equivalence of categories.
\end{enumerate}
%
%
\end{theorem}
\begin{pf} 
The map $\pi:X(D,r)\to X$ is the composition of the maps appearing
in (\ref{tow1}), and the section $\sigma:X\to X(D,r)$ is the composition
of the inclusions in (\ref{inclXj}). The fact that $\mu(r)$ induces an
equivalence of categories follows from Lemma \ref{mainthmj}.
\end{pf}

\section{Cohomological questions}
\label{sec:cohoques}

In this section we address two different questions. First, that of relating
the topology of $W\in\VeqG(X(D,\ur))$ to that of 
$M(W)=(V,\VVV_1,\dots,\VVV_s)\in\PPP(X,D,\ur)$. 
In particular, we obtain a formula for the
first Chern class of $W$ in terms of the first Chern class of $V$ and the
ranks of the elements in the filtrations $\VVV_u$.

Secondly, we compute the Kaehler cone of $X(D,\ur)$ in terms of the
Kaehler cone of $X$ and the classes $[D_1],\dots,[D_s]\in H_2(X)$.
We reduce the computation to the case $D$ smooth and $r=1$, i.e., 
to the problem of relating the Kaehler cone of the blow up of 
$X\times\PP^1$ along $D\times\{[0:1]\}$ to the Kaehler cone
of $X$ and the class $[D]\in H_2(X)$.

\subsection{Deformations near the divisor}
Let $X,X'$ be manifolds, and let $D\subset X$ and
$D'\subset X'$ be smooth divisors.
Let $(V,V_1)\in\PPP(X,D)$ and $(V',V'_1)\in\PPP(X',D')$. We will say that
$(X,D,V,V_1)$ and $(X',D',V',V'_1)$ are isomorphic near the
divisor, and we will write
$$(X,D,V,V_1)\IND(X',D',V',V'_1),$$
if there exist a neighbourhood $U$ (resp. $U'$) of $D$ 
resp. $D'$), a biholomorphism $\phi:U\to U'$ which identifies $D$ with
$D'$, and an isomorphism $\psi:V|_{U}\to\phi^*V'|_{U'}$ which
identifies $V_1$ with $\phi^*V'_1$. If $p:X_D\to X$ and $p':X'_{D'}\to X'$
denote the projections, it follows from Theorem \ref{mainthm} that,
for any pair of objects
$W\in\Veq(X_D)$ and $W'\in\Veq(X'_{D'})$ such that
$\mu(W)\simeq(V,V_1)$ and $\mu(W')\simeq(V',V'_1)$,
$\psi$ induces an isomorphism between 
$W|_{p^{-1}U}$ and $W'|_{{p'}^{-1}U'}$ respectively.

We will say that $(X,D,V,V_1)$ and $(X',D',V',V'_1)$ are directly deformation
equivalent near the divisor, and we will write
$$(X,D,V,V_1)\DDEND(X',D',V',V'_1),$$
if there exists a submersion $q:Y\to B$,
where $B$ is smooth and connected, a
divisor $\Delta\subset Y$ such that $q|_{\Delta}:\Delta\to B$ 
is also a submersion,
a pair $(Z,Z_1)\in\PPP(Y,\Delta)$, and two points $b,b'\in B$
such that 
$$(q^{-1}b,(q|_{\Delta})^{-1}b,Z|_{q^{-1}b},
Z_1|_{(q|_{\Delta})^{-1}b})\IND(X,D,V,V_1)$$ 
and $$(q^{-1}b',(q|_{\Delta})^{-1}b',Z|_{q^{-1}b'},
Z_1|_{(q|_{\Delta})^{-1}b'})\IND(X',D',V',V'_1).$$
We will call deformation equivalence near the divisor, and denote it
by $\DEND$, the equivalence relation induced by $\DDEND$. 

\begin{lemma}
If $(X,D,V,V_1)\DEND(X',D',V',V'_1)$
and we have $(V,V_1)=\mu(W)$, $(V',V'_1)=\mu(W')$ for some
$W\in\Veq(X_D)$ and $W'\in\Veq(X'_{D'})$, then 
there is a $C^{\infty}$ isomorphism between the restrictions of
$W$ and $W'$ to small neighbourhoods
of $p^{-1}D$ and ${p'}^{-1}D'$ respectively.
\label{lemma60}
\end{lemma}
\begin{pf}
Combine Theorem \ref{mainthm} with Ehresmann's Theorem.
\end{pf}

Let $\pi:L\to X$ be the line bundle obtained from $D$ and let
$\sigma_0:X\to L$ be the zero section. Let $\sigma\in H^0(L)$ be 
a section transverse to $\sigma_0$ and such that $\sigma^{-1}(0)=D$.
Let $\pi_D:L|_D\to X$ be the restriction of $\pi$. 
The following result will be useful to 
relate the topology of an object $W\in\Veq(X_D)$ and its image
$(V,V_1)$ by $\mu$.

\begin{lemma}
Let $(V,V_1)\in\PPP(X,D)$ be any object, and let $V_0=V|_D/V_1$.  
Then $$(X,D,V,V_1)\DEND
(L|_D,\sigma_0(D),\pi_D^*(V_0\oplus V_1),\pi_D^*V_1|_{\sigma_0(D)}).$$  
\label{deformacio}
\end{lemma}
\begin{pf}
The proof is split in two steps. We first prove that 
$$(X,D,V,V_1)\DDEND
(L|_D,\sigma_0(D),\pi_D^*V,\pi_D^*V_1|_{\sigma_0(D)}).$$
We use for that the trick of deformation to the normal cone.
Let $B=B(0,2)\subset\CC$, and let $B^{\times}=B\cap\Gm$. 
Let $Y^{\times}$ be the graph of the map $\Sigma:X\times B^{\times}\to
L$ which sends $(x,t)$ to $t^{-1}\sigma(x)$. Let $Y$ be the closure of
$Y^{\times}$ inside $X\times B\times L$, and let $q:=\pi_2:Y\to B$ be the
projection to the second factor. Let ${\Delta}\subset Y$ be the divisor
$\{(y,t,\sigma_0(y))\mid y\in D,\ t\in B\}$. Let us prove that both
$q$ and its restriction to $\Delta$ are submersions. This is clearly
true in $q^{-1}B^{\times}=Y^{\times}$. As for the points in
$q^{-1}(0)$, observe that $\pi_1 q^{-1}(0)\subset D$, where
$\pi_1:X\times B\times L\to X$ is the projection. It thus suffices to
take for any $x\in D$ any neighbourhood $x\in E\subset X$ and study
$\pi_1^{-1}(E)\cap q^{-1}(0)$. So suppose that $E$ is small enough
so that there are coordinates $x_1,\dots,x_n$ on $E$ such that $x$
corresponds to $(0,\dots,0)$ and $D=\{x_1=0\}$. Then we can trivialize
$L|_E\simeq\uC$ in such a way that the section $\sigma$ maps
$(x_1,\dots,x_n)$ to $x_1$. Identifying $X$ with $E$ it follows that
$Y^{\times}=\{(x_1,\dots,x_n,t,t^{-1}x_1)\mid t\in B^{\times}\}$ (in
the last term we just write the fibrewise component of the points in
$L|_E$, using the trivialisation). So $Y=\{(x_1,\dots,x_n,t,y)\mid
t\in B, ty=x_1\}$. But it is then clear that the map
$Y\ni(x_1,\dots,x_n,t,y)\mapsto t\in B$ is a submersion.

Returning to the case of general $X$, define $Z:=\pi_1^*V\to Y$, and
$Z_1:=(\pi_1|_{\Delta})^*V'$. Then it is easy to check that
$$(X,D,V,V_1)\IND(q^{-1}(1),q^{-1}(1)\cap
{\Delta},Z|_{q^{-1}(1)},Z_1|_{q^{-1}(1)\cap {\Delta}})$$ and that
$$(L|_D,\sigma_0(D),\pi_D^*V,\pi_D^*V_1|_{\sigma_0(D)})\IND  
(q^{-1}(0),q^{-1}(0)\cap {\Delta},Z|_{q^{-1}(0)},Z_1|_{q^{-1}(0)\cap
  {\Delta}}).$$

In the second step we prove that
$$(L|_D,\sigma_0(D),\pi_D^*V,\pi_D^*V_1|_{\sigma_0(D)})\DDEND
(L|_D,\sigma_0(D),\pi_D^*(V_0\oplus V_1),\pi_D^*V'|_{\sigma_0(D)}).$$  
This will be a consequence of the following general fact: if
$V_1\subset V$ are vector bundles on $X$ and we define $V_0=V/V_1$, 
then there is a bundle
$Z\to Y=X\times\CC$ such that $Z|_{X\times\{1\}}\simeq V$ and
$Z|_{X\times\{0\}}\simeq V_0\oplus V_1$. Indeed, let us take a 
$C^{\infty}$ splitting $V\simeq V_0\oplus V_1$,
so that the $\ov{\partial}$-operator of $V$ may be written as
$$\ov{\partial}_V=\left(\begin{array}{cc}
\ov{\partial}_{V_0} & \beta \\
0 & \ov{\partial}_{V_1}\end{array}\right),$$
where $\beta\in\Omega^{0,1}(X;V_1\otimes V_0^*)$. Then define
$Z:=\pi_X^*(V_0\oplus V_1)$ with the following
$\ov{\partial}$-operator: 
$$\ov{\partial}_{Z}|_{X\times\{t\}}=\left(\begin{array}{cc}
\ov{\partial}_{V'} & t\beta \\
0 & \ov{\partial}_{V''}\end{array}\right).$$
This clearly satisfies the desired properties.
\end{pf}

\subsection{The functor $\mu$ and topology of vector bundles}
Here we will assume that $X$ is connected.
All the cohomology groups appearing in this section will be taken 
with coefficients in $\RR$.

Let $Y=X\times\PP^1$, and consider the action of $\Gm$ on $Y$ given by
$\theta\cdot(z,[y:w])=(z,[y:\theta w])$. 

\begin{lemma}
There is a $\Gm$-equivariant map $q:X_D\to Y$ which identifies $X_D$
with the blow up of $Y$ along $Z=D\times\{[0:1]\}$ and such that the
exceptional divisor $E=q^{-1}Z$ corresponds to $\Delta^+$.
\label{esunblowup}
\end{lemma}
\begin{pf}
Recall that by definition $X_D\subset\PP(L\oplus\uC\oplus\uC)$ and
that $p:X_D\to X$ denotes the projection.
Let us identify $\PP(\{0\}\oplus\uC\oplus\uC)$ with
$\PP(\uC\oplus\uC)=X\times\PP^1=Y$ by 
mapping $[0:y:w]\in\PP(\{0\}\oplus\uC\oplus\uC)_z$ to $(z,[y:w])$.
This bijection is $\Gm$-equivariant. With this in mind we will define
the desired map $q$ from $X_D$ to $\PP(\{0\}\oplus\uC\oplus\uC)$. 

Let $z\in X$ and let $v\in p^{-1}(z)$ be any point. Let
$u=p^{-1}(z)\cap X_2$, and let
$\Lambda_v\subset\PP(L\oplus\uC\oplus\uC)_z$ be the line passing
through $u$ and $v$ if $u\neq v$, or the tangent of $p^{-1}(z)$ at $u$
if $u=v$.  Then we set 
$$q(v):=\PP(\{0\}\oplus\uC\oplus\uC)_z\cap\Lambda_v.$$
Let us check that $q$ is the blow up of $Y$ along $Z$ and that $q$ is
$\Gm$-equivariant. For this it
suffices to work locally in $X$, so that we may assume that $X=\CC^n$
and $D=\{0\}\times\CC^{n-1}$. Then $X_D=H_n$, and the map $q$ can be
described as follows. A point $v=[x:y:w]\in \{z\}\times\PP(\CC^3)$
satisfying $xy=w^2t_1$, where $z=(t_1,\dots,t_n)$, is mapped to
$$q(v)=\left\{
\begin{array}{ll}
(z,[y:w]) & \text{ if $[x:y:w]\neq[1:0:0]$ } \\
(z,[0:1]) & \text{ if $[x:y:w]=[1:0:0]$ }
\end{array}\right.$$
From this it is clear that $q$ is $\Gm$-equivariant.
It also follows that for any $(z,[y:w])$ such that $y\neq 0$
the preimage $r^{-1}(z,[0:y:w])$ is equal to
$\{(z,[y^{-1}w^2t_1:y:w])\}$, so that the restriction of $q$ to
$q^{-1}\{y\neq 0\}$ is an isomorphism. 
The blow up of $\{w\neq 0\}$ along $Z$ (note that $Z\subset\{w\neq
0\}$) is by definition
$$B=\{(t_1,\dots,t_n,[y:1],[a:b])\mid at_1=by\}.$$
It now suffices to construct a map $\Psi:B\to X_D$ which induces a
biholomorphism between $B$ and $\Psi(B)$ and which makes the following
diagram commute:
$$\xymatrix{
B\ar[rr]^{\Psi}\ar[rd]_{r} && X_D\ar[ld]^{q} \\
& Y,}$$ 
where $r$ is the composition of the blow up map $B\to\{w\neq 0\}$ with
the inclusion $\{w\neq 0\}\subset Y$. The following definition suits
our needs:
$$\Psi(t_1,\dots,t_n,[y:1],[a:b])=\left\{
\begin{array}{ll}
(t_1,\dots,t_n,[t_1:y^2:y]) & \text{ if $(t_1,y)\neq (0,0),$}\\
(t_1,\dots,t_n,[b:0:a]) & \text{ if $(t_1,y)=(0,0).$}
\end{array}\right.$$
From this explicit computation and (\ref{descriuD}) it follows that
the exceptional divisor of $q$ is $\Delta^+$.  
\end{pf}
Let $N$ be the normal bundle of the inclusion $Z\subset Y$.
Then $E=\PP(N)$.

\begin{lemma}
There is an isomorphism of cohomology with real coefficients
$$H^2(X_D)=H^2(X)\oplus H^2(\PP^1)\oplus \RR\la t\ra,$$
where $t$ is the first Chern class the line bundle associated to $E$.
\label{H2XD}
\end{lemma}
\begin{pf}
This follows from Lemma \ref{esunblowup} and a standard application of
Mayer--Vietoris.
\end{pf}

Denote by $i:E\to X_D$ the inclusion of the exceptional
divisor, and let $N_E\to E$ be the normal bundle. 
By a slight abuse of notation, if $a\in H^*(E)$, $t\cup a$ will
denote the image of $a$ by the composition
\begin{equation}
H^*(E)\to H^*(N_E,N_E\setminus\{0\})\to H^*(X_D,X_D\setminus E)
\to H^*(X_D),
\label{cuprar}
\end{equation}
where the first map is Thom's isomorphism, the second one is
excission, and the third one is the map induced by the inclusion
$(X_D,\emptyset)\subset (X_D,X_D\setminus E)$. Similarly, for any
natural number $n\geq 1$, $t^n\cup a$ will sometimes denote
$t\cup(i^*t^{n-1}\cup a)$, where the first $\cup$ is the one 
defined by (\ref{cuprar}) and the second one is the usual one in
$H^*(E)$. For example, if $a=i^*b$ for some $b\in H^*(X_D)$, then
\begin{equation}
t^n\cup a=t^n\cup b, 
\label{formulacuprar}
\end{equation}
where the $\cup$ in the RHS is the usual one
in $H^*(X_D)$. Finally, we extend the definition by linearity in order
to give a sense to expressions of the form $P(t)\cup a$, where $P$ is
a polynomial satisfying $P(0)=0$.


\begin{lemma}
Let $W\in\Veq(X_D)$, and let $(V,V_1):=\mu(W)$. The Chern character of $W$ is
$$\ch W = p^* \ch V + (e^t-1)\cup p^*\ch V_1.$$
\label{ChernW}
\end{lemma}
\begin{pf}
Take some $W\in\Veq(X_D)$ and set $(V,V_1):=\mu(W)$. 
Let also $W_0\in\Veq(X_D)$ be such that $\mu(W_0)\simeq (V,0)$.
By Corollary \ref{peszerotrivial}, $W_0$ is isomorphic to $p^*V$ with the
trivial action of $\Gm$. So 
\begin{equation}
\ch W_0=p^* \ch V.
\label{chW0}
\end{equation}
Let $X_D^*=X_D\setminus p^{-1}D$. Consider the cohomology long exact
sequence for the pair $(X_D,X_D^*)$:
$$\dots\longrightarrow H^k(X_D,X_D^*)\stackrel{i}{\longrightarrow} 
H^k(X_D)\stackrel{j}{\longrightarrow} H^k(X_D^*)\longrightarrow\dots$$
By Corollary \ref{peszerotrivial} the
restrictions $W|_{X_D^*}$ and $W_0|_{X_D^*}$ are isomorphic. Hence 
$j(\ch W-\ch W_0)=0$, which implies that $\ch W-\ch W_0=i(d_{W,W_0})$ 
for some class $d_{W,W_0}\in H^*(X_D,X_D^*)$.
On the other hand, if $B$ is any neighbourhood of $D$, we have an
isomorphism induced by excission
$$H^*(X_D,X_D^*)\simeq H^*(p^{-1}B,p^{-1}(B\setminus D)),$$
so that to compute $\ch W-\ch W_0$ it suffices to compute
$\ch W|_{p^{-1}B}-\ch W_0|_{p^{-1}B}$ for any neighbourhood $B$ of
$D$. 

Now, by Lemmae \ref{lemma60} and
\ref{deformacio} we can reduce the computation to
the case in which $V$ splits as $V_a\oplus V_b$ and $V_1$ is equal to
the restriction $V_a|_D$. Then, for any
$W_a,W_b\in\Veq(X_D)$ such that $\mu(W_a)$ (resp. $\mu(W_b)$)
is isomorphic to $(V_a,V_a|_D)$ (resp. $(V_b,0)$)
there is a $C^{\infty}$ isomorphism 
$W|_B=W_a\oplus W_b$. By Corollary \ref{peszerotrivial}
we have $\ch W_b=p^*\ch V_b$. Let $\LLL\to X_D$ be the line bundle
constructed in Subsection \ref{secondthing}. One can check that 
$W_a'=W_a\otimes\LLL|_{p^{-1}B}^{-1}$ belongs to $\Veq(p^{-1}B)$ 
and that all its weights are zero. 
Hence, by Corollary \ref{peszerotrivial}, we have isomorphisms
$$W_a'\simeq p^*W_a'|_{X_2\cap p^{-1}B}\simeq p^*W_a|_{X_2\cap p^{-1}B}$$ 
(the last isomorphism follows
from the existence of the nowhere zero section $\psi\in
H^0(X_2;\LLL|_{X_2})$). So we can compute the Chern character:
$\ch W_a=\ch\LLL\cup p^*\ch V_a$. But by construction the first Chern
class of $\LLL$ is $t$, so  $\ch W_a=e^t\cup p^*\ch V_a$ and
consequently $\ch W|_{p^{-1}B}=p^*\ch V_b+e^t\cup p^*\ch V_a$. Hence, 
\begin{align}
\ch W|_{p^{-1}B}-\ch W_0|_{p^{-1}B} &=
p^*\ch V_b+e^t\cup p^*\ch V_a - p^*\ch(V_a\oplus V_b) \notag \\
&= (e^t-1)\cup p^* \ch V_a \notag \\
&= (e^t-1)\cup p^* \ch V_1, \label{laresta}
\end{align}
where the last inequality follows from (\ref{formulacuprar})
and in the last line the cup product is the one defined by (\ref{cuprar}).
Combining (\ref{chW0}) with (\ref{laresta}) the result follows.
\end{pf}

\begin{corollary}
Let $W\in\Veq(X,D)$, and let $(V,V_1)=\mu(W)$.
Then $$c_1(W)=p^* c_1(V)+\rk V_1t.$$
\label{c1W}
\end{corollary}

\subsection{The Kaehler cone of certain blow ups}

\subsubsection{}
\label{constrcon}
In order to compute the Kaehler cone of $X_D$ we will prove a general
result describing the Kaehler cone of the blow up of a product of
compact Kaehler manifolds $X\times X'$ along $D\times D'$, where
$D\subset X$ and $D'\subset X'$ are smooth divisors.
This is done in this subsection, the main result being 
stated in Theorem \ref{KXX}.

Let $N$ be a manifold and $\pi:\Lambda\to N$ a Hermitian line bundle, and
assume that a compact connected Lie group $T$ acts on $\Lambda$ on the left,
linearly on the fibres and respecting the metric. Let $\tlie=\Lie T$.
For any $s\in\tlie$ denote by $\fX_\Lambda(s)$ the vector field
generated by the infinitesimal action of $s$ on $\Lambda$ and by
$\fX_N(s)=d\pi\fX_\Lambda(s)$ the vector field generated on $N$. 

Fix a $T$-invariant Hermitian connection $B$ on $\Lambda$. This defines for
any vector field $\fX\in\Gamma(TN)$ a lift
$\sigma_B(\fX)\in\Gamma(T\Lambda)$. We define the moment of the action of
$T$ on $\Lambda$ w.r.t. $B$ to be the map
$\Omega_B:\tlie\to\Omega^0(N;\imag\RR)$ which 
assigns to $s\in\tlie$ the field
$$\Omega_B(s):=\sigma_B(\fX_N(s))-\fX_\Lambda(s).$$
This is easily seen to be a vertical vector field and, using the
canonical identification $T\Lambda_v\simeq\pi^*\Lambda$, it is in fact
a Hermitian 
linear map, so we can view $\Omega_B(s)\in\Omega^0(N;\imag\RR)$.
A straightforward computation proves that for any field $\fX\in\Gamma(TN)$
we have
\begin{equation}
d\Omega_B(s)=i(\fX_N(s))F_B,
\label{primeraest}
\end{equation}
where $F_B\in\Omega^2(N;\imag\RR)$ is the curvature of $B$. 

Now assume that $M$ is another manifold, and that $P\to M$ is a
$T$-principal bundle. We will construct a connection $\nabla$ on
the line bundle
$$\pi_M:P\times_T \Lambda\to P\times_T N,$$
and for that it suffices to specify a suitable retraction
$\nu_\nabla:T(P\times_T \Lambda)\to\Ker d\pi_M$. 
Let $p:P\times_T \Lambda\to M$ be the projection. Take a connection $A$ on
$P$. Then we get a retraction $\nu_A:T(P\times_T \Lambda)\to \Ker
dp=P\times_T T\Lambda$. On the other hand, the connection $B$ gives a
$T$-equivariant retraction $\nu_B:T\Lambda\to \Ker d\pi$ which consequently
extends to a retraction 
$$\nu_B':P\times_T T\Lambda\to P\times_T\Ker d\pi=\Ker d\pi_M.$$
We then define
$$\nu_\nabla:=\nu_B'\circ\nu_A.$$
It is easy to see that this defines a connection $\nabla$ on
$P\times_T \Lambda\to P\times_T N$.

\begin{lemma}
The curvature of $\nabla$ is 
$$F_\nabla={\nu_A'}^*F_B-\Omega_B(q^*F_A),$$
where $q:P\times_T N\to M$ is the projection and $\nu_A':T(P\times_T
N)\to \Ker dq$ is the retraction induced by $A$.
\label{lacurvatura}
\end{lemma}
\begin{pf}
This follows from an easy local computation (see for example Lemma
A.1.4 in \cite{M}). Note, by the way, that $F_{\nabla}$ is the
coupling form of the bundle $P\times_TN$, see \cite{GLS}.
\end{pf}

\subsubsection{}
Let $N=\PP^1$ with the action of $T=S^1\times S^1$ given by
\begin{equation}
(\theta_1,\theta_2)[x:y]=[\theta_1x:\theta_2y].
\label{casP1}
\end{equation}
This action lifts to $\Lambda:=\OOO(-1)\to\PP^1$ by defining
$(\theta_1,\theta_2)(\lambda x,\lambda y)=(\theta_1\lambda
x,\theta_2\lambda y)$, where we identify the fibre of $\OOO(-1)$ over
$[x:y]$ with $\{(\lambda x,\lambda y)\mid \lambda\in\CC\}$. Take the
$T$-invariant metric $|(\lambda x,\lambda y)|^2=|\lambda x|^2+|\lambda
y|^2$ on $\Lambda$, and let $B$ be its Chern connection. Then $B$ is
$T$-invariant, and its curvature satisfies
$$F_B=2\pi\imag \omega_{\PP^1},$$
where $\omega_{\PP^1}$ is the Fubini-Study symplectic form on $\PP^1$.
Identifying $\tlie=\Lie T\simeq\imag\RR\oplus\imag\RR$, let
$t_1=(i,0)\in\tlie$, $t_2=(0,i)\in\tlie$, 
and let $\tau_1,\tau_2$ be the dual basis of $t_1,t_2$. Define
$$\mu([x:y]):=\frac{\tau_1|x|^2+\tau_2|y|^2}{|x|^2+|y|^2}\in\tlie^*.$$ 
Then $\mu$ is a symplectic moment map for the action of $T$ on $\PP^1$
with respect to $\omega_{\PP^1}$.
It follows from (\ref{primeraest}) that
$\Omega_B=\mu+C$ for some constant $C\in\tlie^*$. 
Let us compute this constant.
The subgroup $G=S^1\times\{1\}\subset T$ acts trivially on the fibre
$\Lambda_{[0:1]}$ and the action on the tangent space
$T_{[0:1]}\PP^1$ has weight $1$. So for any $x\in\Lambda_{[0:1]}$ the
tangent space $T_x\Lambda$ splits by the action of $G$ in weights
$0$ and $1$. And, since the connection $B$ is $G$-invariant, the
horizontal lift of $T_{[0:1]}\PP^1$ must be the inclusion in the
summand of weight $1$ in $T_x\Lambda$. But the value of the vector
field $\fX_{\Lambda}(t_1)$ at $x$ lies also inside the same
summand. Consequently, $\Omega(t_1)([0:1])=0$. Similarly one checks
that $\Omega(t_2)([1:0])=0$, so we deduce that $C=0$, i.e., 
$\Omega_B=\mu.$

\subsubsection{}
For any manifold $Z$ we will denote by $K(Z)\subset
H^{1,1}(Z;\CC)\cap H^2(Z;\RR)$ the Kaehler cone of $Z$.

\begin{theorem}
Let $X,X'$ be Kaehler manifolds, with $X'$ simply connected, and let
$D\subset X$, $D'\subset X'$ be smooth divisors. Let $c=c_1(D)$ and
$c'=c_1(D')$. Let $Y$ be the blow up of $X\times X'$ at $D\times
D'$. Then we have 
\begin{equation}
H^2(Y;\RR)=H^2(X;\RR)\oplus H^2(X';\RR)\oplus \RR\la t\ra,
\label{cincest}
\end{equation}
where $t$ is the first Chern class of the exceptional divisor. The
Kaehler cone of $Y$ is
$$K(Y)=\left\{(w,w',bt)\in H^2(Y;\RR)\left|
\begin{array}{l}
w\in K(X),\ w+bc\in K(X), \\
w'\in K(X'),\ w'+bc'\in K(X'), \\
\text{and }b<0.
\end{array}
\right.\right\}.$$
Suppose now that $S^1$ acts on $X'$ keeping $D'$ fixed, so that there
is an induced action of $S^1$ on $X\times X'$ and the blow up
$Y$. Then any class in $K(Y)$ can be represented by a $S^1$-invariant
Kaehler form.
\label{KXX}
\end{theorem}
\begin{pf}
Formula (\ref{cincest}) follows from K{\"u}nneth (using
$H^1(X';\RR)=0$) and 
an easy Mayer--Vietoris argument, exactly as in Lemma \ref{H2XD}. 
Let $L\to X$ and $L'\to
X'$ be the line bundles defined by $D$ and $D'$. Since the normal
bundle 
$N_{D\times D'\mid X\times X'}=L|_D\boxplus L'|_{D'},$
there is an embedding $$i_Y:Y\to\PP(L\boxplus L')=:Q.$$
More precisely, let $\sigma\in H^0(L)$ and $\sigma'\in H^0(L')$ be
sections transverse to zero such that $\sigma^{-1}(0)=D$ and
${\sigma'}^{-1}(0)=D'$. Then $i_Y(Y)$ coincides with the closure in
$Q$ of the image of the section
$$\begin{array}{rcl}
X\times X'\setminus D\times D' & \longrightarrow & Q \\
(x,x') & \mapsto & [\sigma(x):\sigma'(x')].
\end{array}$$

Let $x_0\in D$ and $x_0'\in D'$. Let $i_0:X\setminus D\to
Q|_{X\times\{x_0'\}}$ and $i_0':X'\setminus D'\to
Q|_{\{x_0\}\times X'}$ be the maps
$i_0(x)=[\sigma(x):0]$ and $i_0'(x')=[0:\sigma'(x')]$. Both $i_0$ and
$i_0'$ extend to give embeddings $i:X\to Y$ and $i':X'\to Y'$
satisfying
$$i^*(w,w',bt)=w+bc\qquad\text{ and }\qquad {i'}^*(w,w',bt)=w'+bc'.$$
Taking now $x_0\in X\setminus D$ and $x_0'\in X'\setminus D'$ and
repeating the same construction we get embeddings $j:X\to Y$ and
$j':X'\to Y$ such that $\Im i_Yj\subset Q|_{X\times\{x_0'\}}$ and
$\Im i_Yj'\subset Q|_{\{x_0\}\times X'}$ and satisfying
$$j^*(w,w',bt)=w\qquad\text{ and }\qquad {j'}^*(w,w',bt)=w'.$$

Now suppose that $(w,w',bt)\in K(Y)$. Then, since $i,i',j,j'$ are
embeddings, $w\in K(X)$, $w+bc\in K(X)$, $w'\in K(X')$ and $w'+bc'\in
K(X')$ (indeed, the pullback by an embedding of a Kaehler form is a
Kaehler form). 
Let $\pi:Y\to X\times X'$ denote the projection. To check that $b<0$,
take any point $z=(x,x')\in D\times D'$. Then
$Z:=\pi^{-1}(z)=\PP(L_x\oplus L'_{x'})$ and the restriction of $t$ to
$Z$ is $t|_Z=c_1(\OOO_{\PP(L_x\oplus L'_{x'})}(-1))$. Hence,
$$\int_Z w+w'+bt=\int_Z bt=-b,$$
and this must be positive if $(w,w',bt)$ is a Kaehler class.

Conversely, let $(w,w',bt)\in H^2(Y;\RR)$ satisfy
$w\in K(X)$, $w+bc\in K(X)$, $w'\in K(X')$, $w'+bc'\in
K(X')$ and $b<0$. We will construct a Kaehler form
$\kappa\in\Omega^2(Y)$ representing the class $w+w'+bt$.
Let $\omega_0,\omega_1\in\Omega^2(X)$ and
$\omega'_0,\omega'_1\in\Omega^2(X')$ be Kaehler forms satisfying
$[\omega_0]=w$, $[\omega_1]=w+bc$, $[\omega'_0]=w'$ and
$[\omega'_1]=w'+bc'$. Define $F=-2\pi\imag b^{-1}(\omega_1-\omega_0)$ and
$F'=-2\pi\imag b^{-1}(\omega'_1-\omega'_0)$. 

Fix metrics on $L$ and $L'$. 
One can then take Hermitian connections $C$ on
$L$ and $C'$ on $L'$ such that the curvatures satisfy
$$F_C=F\qquad\text{ and }\qquad F_{C'}=F'.$$
Let $P_L\subset L$ and $P_{L'}\subset L'$ be the unit length
vectors. Both $P_L$ and $P_{L'}$ are $S^1$ principal bundles.
Define $M:=X\times X'$ and $P=\pi_X^*P_L\times \pi_{X'}^*P_{L'}$. Then
$P$ is a $T:=S^1\times S^1$ principal bundle. 

Using the action of $T$ on $\PP^1$ defined by (\ref{casP1}), we have
$P\times_T\PP^1=Q$. By the K{\"u}nneth and Leray--Hirsch theorems and
the fact that $H^1(\PP^1;\RR)=0$ we have 
$$H^2(Q;\RR)=H^2(X;\RR)\oplus H^2(X';\RR)\oplus\RR\la t_Q\ra,$$
where $t_Q=c_1(\OOO_{Q}(-1))$ is the relative canonical bundle. 
Taking the bundle $\Lambda=\OOO_{\PP^1}(-1)$ 
and the lift of the action of $T$ defined
in \ref{casP1}, we also have
$$\OOO_{Q}(-1)=P\times_T\Lambda.$$
Let $A$ be the connection on $P$ induced by the pullbacks of the
connections $C$ and $C'$. Then
$F_A=F+F'\in\Omega^2(X\times X')$ (we omit the
pullbacks). Taking on $\Lambda$ the connection specified in
\ref{casP1}, we can now apply the construction in \ref{constrcon} to get
a connection $\nabla$ on $\OOO_{Q}(-1)$.
Then we define 
$$\kappa_Q:=\frac{\imag}{2\pi}(\pi^*(\omega_0+\omega_0')+bF_{\nabla}).$$
It is clear that this form represents $[\kappa_Q]=(w,w',bt_Q)$. 
Let us check that it is a Kaehler form, i.e., that it is positive.
Take a point $y\in Q$ and a nonzero $0\neq v\in T_yQ$. We want to
prove that $$\kappa_Q(v,iv)>0.$$ Denote by $\pi_Q:Q\to X\times X'$ the
projection, and write $\pi_Q(y)=(x,x')$ and $d\pi_Q(v)=(u,u')\in
T_xX\oplus T_{x'}X'$. Let also $v_0=\nu_A'(v)$ be the vertical
projection of $v$. Suppose that $y=[y_0:y_1]\in\PP(L_x\oplus
L'_{x'})$. By Lemma \ref{lacurvatura} and (\ref{cincest}) we have
\begin{align*}
\kappa_Q(v,iv)&=\frac{\imag}{2\pi}
\left(bF_B(v_0,iv_0)+b\frac{|y_0|^2}{|y_0|^2+|y_1|^2} F(u,iu)
+b\frac{|y_1|^2}{|y_0|^2+|y_1|^2} F'(u',iu')\right) \\
&+\omega_0(u,iu)+\omega_0'(u',iu').
\end{align*}
By construction, $\frac{\imag}{2\pi}bF_B(v_0,iv_0)>0$. As for the
other terms, recall that 
$$\frac{\imag b}{2\pi}F=\omega_1-\omega_0
\qquad\text{ and }\qquad
\frac{\imag b}{2\pi}F'=\omega'_1-\omega'_0,$$
so the remaining summands in $\kappa_Q(v,iv)$ can be written as
$$(\lambda\omega_0+(1-\lambda)\omega_1)(u,iu)+
((1-\lambda)\omega'_0+\lambda\omega'_1)(u',iu'),$$  
for some $0\leq \lambda\leq 1$, and this is positive since by assumption
$\omega_0,\omega_1,\omega'_0,\omega'_1$ are Kaehler forms.

Finally, we set $\kappa:=i_Y^*\kappa_Q$. Since $i_Y$ is an embedding,
$\kappa$ is a Kaehler form, and since $i_Y^*t_Q=t$, it follows that
$\kappa$ represents the class $(w,w',bt)$ as desired.

It remains to prove the last statement on existence of $S^1$-invariant
Kaehler forms. This follows from the standard averaging trick.
\end{pf}


\subsection{The Kaehler cone of $X(D,r)$}
\label{KXDr}

In the rest of the present section we will assume that $X$ is Kaehler.
Recall that $X(D,r)$ fits in a tower of maps
$$X(D,r)=Y_r\to Y_{r-1}\to\dots\to Y_0=X.$$
We will describe the Kaehler cone of $Y_j$ for any $1\leq j\leq r$.
We first define (using induction on $j$) cohomology
classes $$d_{j,1},\dots,d_{j,j},t_{j,1},\dots,t_{j,j}$$ in $H^2(Y_j)$.
Recall that $Y_1=X_D$ is isomorphic to the blow up of
$Y=X\times\PP^1$ along $D\times\{[0:1]\}$
(see Lemma \ref{esunblowup}). Let $q:X_D\to Y$ be the blow up
map, and let 
$\delta=1\otimes PD[\PP^1]\in H^0(X)\otimes H^2(\PP^1)\subset H^2(Y)$.
Then we set $d_{1,1}:=q^*\delta$, and we define $t_{1,1}$ to be first
Chern class of the exceptional divisor 
$\Delta^+_1=q^{-1}(D\times\{[0:1]\})\subset Y_1$.
Now assume that $1\leq j<r$ and that the classes
$d_{j,1},\dots,d_{j,j},t_{j,1},\dots,t_{j,j}\in H^2(Y_j)$ have already
been defined. Recall that $Y_{j+1}=(Y_j)_{\Delta^+_j}$. Using again
Lemma \ref{esunblowup} we can factor the projection map
$p_{j+1}:Y_{j+1}\to Y_j$ as follows
$$\xymatrix{Y_{j+1} \ar[rd]^{q_{j+1}} \ar[dd]_{p_{j+1}}\\
& Y_j\times\PP^1 \ar[ld]^{o_j} \\
Y_j,}$$
where $q_{j+1}$ is the blow up map of $Y_j\times\PP^1$ along
$\Delta^+_j\times\{[0:1]\}$. We then define for any $1\leq i\leq j$
$$d_{j+1,i}:=p_{j+1}^*d_{j,i}\qquad\text{ and }\qquad
t_{j+1,i}:=p_{j+1}^*t_{j,i}.$$
We also set $d_{j+1,j+1}:=q_{j+1}^*(1\otimes PD[\PP^1])$ and
$t_{j+1,j+1}=c_1(\Delta^+_{j+1})$.
In the sequel, by an abuse of notation, we will omit the first subindex
in the variables $d_{j,i},t_{j,i}$, so that $d_i$ and $t_i$
will denote cohomology classes in $Y_j$ for every $j\geq i$).

For any $1\leq j\leq r$ let $\rho_j$ denote the composition
$p_{j}\circ\dots\circ p_1:Y_j\to X$. 

\begin{lemma}
Let $1\leq j\leq r$. The map $\rho_j^*:H^2(X)\to H^2(Y_j)$ is injective
and, identifying $H^2(X)$ with its image in $H^2(Y_j)$ by this map, we have
$$H^2(Y_j)=H^2(X)\oplus\RR\la d_1,\dots,d_j,t_1,\dots,t_j\ra.$$
\label{H2XDj}
\end{lemma}
\begin{pf}
It follows from applying recursively Lemma \ref{H2XD}.
\end{pf}

\begin{lemma}
Let $1\leq j\leq r$.
The map $\rho_j^*:H^2(X)\to H^2(Y_j)$ is injective
and, identifying $H^2(X)$ with its image in $H^2(Y_j)$ by this map, we have
$$H^2(Y_j)=H^2(X)\oplus\RR\la d_1,\dots,d_j,t_1,\dots,t_j\ra.$$
Let $c=c_1(D)\in H^2(X)$.
The Kaehler cone of $Y_j$ is equal to
$$K(Y_j)=\left\{(\omega,a,b)\in H^2(Y_j)
\ \left| 
\begin{array}{l}
a\in\RR^j,\ b\in\RR^j,\\
\omega\in K(X),
\end{array}
\begin{array}{l}
\text{satisfying, for any $1\leq i\leq j$,} \\
a_i>0,\ b_i<0, \\
\omega+(b_1+\dots+b_i)c\in K(X), \\
a_i+b_i+b_{i+1}+\dots+b_j>0.
\end{array} 
\right.\right\},$$
where for any $\omega\in H^2(X)$, $a=(a_1,\dots,a_j)\in\RR^j$ and
$b=(b_1,\dots,b_j)\in\RR^j$ the triple $(w,a,b)$ denotes 
$$(\omega,a,b):=\omega+\sum_{i=1}^j a_id_i+b_it_i\in H^2(Y_j).$$
\label{KXr}
\end{lemma}
\begin{pf}
Use induction on $j$ and apply at each step Theorem \ref{KXX}
with $X'=\PP^1$.
\end{pf}

\section{Slope of vector bundles (smooth divisor)}
\label{slopesmooth}
The aim of this section is to study the behaviour of the slope of
vector bundles under the functor $\mu(r)$. 
We will throughout assume that $X$ is Kaehler and compact.
We will use the notations
of the preceding sections, specially \ref{KXDr} (recall, by the way,
that the notion of slope is intrinsically related to the Kaehler
cone of the base manifold --- see (\ref{formslope})).
The main result will be stated in Theorem \ref{calculgrau}.
Before, we state and prove some technical lemmae.

Denote by $\partial_l:\Delta^+_l\hookrightarrow Y_l$ the inclusion
of the exceptional divisor, and let $m_l:\Delta^+_l\to\Delta^+_{l-1}$ 
be the restriction of the projection $p_l$. We then have a chain of
maps
\begin{equation}
\Delta^+_r\stackrel{m_r}{\longrightarrow}
\Delta^+_{r-1}\stackrel{m_{r-1}}{\longrightarrow}
\dots\stackrel{m_{1}}{\longrightarrow}\Delta^+_0=D.
\label{torreDelta}
\end{equation}
We now list some properties of cohomology classes defined in 
\ref{calculgrau}.

\begin{lemma}
\label{propdt}
(1) $d_l^2=0$ and $\partial_l^*d_l=0$;

(2) for any $\alpha\in H^*(Y_l)$, 
$\la t_l\alpha,[Y_l]\ra=\la\partial_l^*\alpha,[\Delta^+_l]\ra$;

(3) $t_l^2=-t_lt_{l-1}$ inside $H^*(\Delta^+_l)$;
for any $\alpha\in H^*(\Delta^+_{l-1})$, 
$$\la t_l m_l^*\alpha,[\Delta^+_l]\ra=-\la\alpha,[\Delta^+_{l-1}]\ra,$$
and $\la m_l^*\alpha,[\Delta^+_l]\ra=0$;

(4) suppose that $\alpha\in H^*(Y_l)$ and that
$\la\alpha,[Y_l]\ra\neq 0$ and $\alpha$ is not of the form 
$t_l\alpha'$ for some $\alpha'\in H^*(Y_l)$, then 
$\alpha=d_lp_l^*\alpha''$ for some $\alpha''\in H^*(Y_{l-1})$, and
$\la\alpha,[Y_l]\ra =\la\alpha'',[Y_{l-1}]\ra$.
\end{lemma}
\begin{pf}
(1) The first claim is obvious, and the second one follows from
the commutativity of the diagram
$$\xymatrix{\Delta^+_l\ar[r]^{\partial_l}\ar[d]_{m_l}
& Y_l\ar[d]^{q_l} \\
\Delta^+_{l-1}\times\{[0:1]\}\ar[r]^{n_l} & Y_{l-1}\times\PP^1,}$$
and the fact that $n_l^* PD(1\otimes[\PP^1])=0$.
(2) This holds because $t_l$ is the Poincar\'e dual of the
fundamental class $[\Delta^+_l]$.
(3) Since $\Delta^+_l$ is the exceptional divisor of the blow 
up of $Y_{l-1}\times\PP^1$ along $\Delta^+_{l-1}\times\{[0:1]\}$, we
can identify $\Delta^+_l=\PP(N^+_{l-1}\oplus\uC)$ (recall that
$N^+_{l-1}\to\Delta^+_{l-1}\subset Y_{l-1}$ is the normal bundle), and
the map $m_l$ is the projection  $\PP(N^+_{l-1}\oplus\uC)\to\Delta^+_{l-1}$; 
with this in mind, we apply
Leray's theorem together with the fact that $c_1(N^+_{l-1})=t_{l-1}$.
(4) If $\alpha$ cannot be written $\alpha=t_l\alpha'$, then certainly
$\alpha=q_l^*\beta$ for some $\beta\in H^*(Y_{l-1}\times\PP^1)$;
furthermore, since the map $q_l$ is birational, we have
$\la\alpha,Y_l\ra=\la\beta,Y_{l-1}\times\PP^1\ra$. Now the claim
follows from K\"unneth's theorem.
\end{pf}

For any $I=(i_1,\dots,i_r)\in\ZZ_{\geq 0}^r$ we define
$$|I|=\sum_{l=1}^r i_l\qquad\text{ and }\qquad
\sigma(I)=\sum_{l=1}^r i_l 2^{l-1}.$$
We will use the standard multiindex notation, so that for example
$d^I$ will denote $\prod_{l=1}^r d_{l}^{i_l}$.

\begin{lemma}
\label{propIJ}
Suppose that for some $I=(i_1,\dots,i_r),J=(j_1,\dots,j_r)\in
\ZZ_{\geq 0}^r$ and $\eta\in H^*(X)$
we have $\la d^It^J\eta,[X(D,r)]\ra\neq 0$.
Then there exists some $0\leq l\leq r$ such that

(a) if $l<r$ then $I=(0,\dots,0,1,\dots,1)$ ($1$'s in the last $r-l=|I|$ 
positions and $0$'s everywhere else), and
if $l=r$ then $I=(0,\dots,0)$; 
besides, for any $l<\mu\leq r$, $i_{\mu}=0$;

(b) for any $1\leq \mu\leq l$ we have 
$\sum_{\mu\leq k\leq l}j_k\geq 2+(l-\mu)$;

(c) we have $|I|+|J|\geq r$ and $\sigma(I)+\sigma(J)\geq 2^{r}-1$, 
with equality in any of the two if and only
if $J=0$ and $I=(1,\dots,1)$ (hence, the case $l=r$); in this situation
\begin{equation}
\la d^It^J\eta,[X(D,r)]\ra=\la d^I\eta,[X(D,r)]\ra=
\la\eta,[X]\ra.
\label{formulalemma2}
\end{equation}
\end{lemma}
\begin{pf}
Throughout the proof we will refer to statements (1)-(4) in Proposition
\ref{propdt}.
Applying (1) and (4) from $l=r$ downwards
as many times as possible (i.e., as long as we don't
find some nonzero $i_l$ or we arrive at $l=0$), we deduce
that $I$ ends with a sequence of $1$'s of length $r-l$, where 
$0\leq l\leq r$. Furthermore, the last $r-l$ positions of $J$ vanish.
Now suppose that $l\geq 0$. Then $i_l\neq 0$, and defining 
$I'=(i_1',\dots,i_r')\in\ZZ_{\geq 0}^r$ by $i'_j=i_j-\delta_{lj}$
we deduce from (4) and (2) that
$$\la d^Jt^I\eta,[X(D,r)]\ra=\la t^I\eta,[Y_j]\ra=
\la t^{I'},[\Delta^+_l].$$
Then (1) tells us that the first $l$ positions of $I$ have to vanish,
so we are done with (a). We now apply successively (3) 
to descend step by step the tower (\ref{torreDelta}), until we arrive
at $D$, and we prove (b).
We now prove (c). First of all, it is clear that if $J=0$ and
$I=(1,\dots,1)$ then $\sigma(I)+\sigma(J)=2^{r-1}-1$.
If this is not the case then $1\leq l\leq r$, and $J$ satisfies (b). 
It follows from this that 
$$\sigma(J)=\sum_{1\leq k\leq l}j_k2^{k-1}>2^l.$$
Let us see why.
If $J=J_0:=(1,1,\dots,1,2,0,\dots,0)$ then we are done.
Now, no matter what $J$ is, we can modify it by a sequence of
moves in such a way that (b) is preserved all the time, $\sigma(J)$ does
not increase, and at the end we arrive at $J_0$. 
Indeed, by (c) we have in general $j_l\geq 2$. If $j_l>2$, then we shift 
$j_l-2$ units from position $l$ to position $l-1$. 
This move preserves (c) and makes $\sigma(J)$ decrease.
Next we look at $j_{l-1}$, which is $\geq 1$ by (c). If $j_{l-1}>1$, then
we shift $j_{l-1}-1$ units from position $l-1$ to position $l-2$.
And so on, until $J=(a,1,\dots,1,2,0,\dots,0)$, where $a\geq 1$.
Then we substitute $a$ by $1$, and we are done.
\end{pf}

\begin{lemma}
\label{propIJ2}
Suppose that $I=(i_1,\dots,i_r),J=(j_1,\dots,j_r)\in
\ZZ_{\geq 0}^r$ and $\eta\in H^*(X)$.
If for some $1\leq l\leq r$ we have
$\la t_ld^It^J\eta,[X(D,r)]\ra\neq 0$ then
$$|I|+|J|\geq r\qquad\text{and}\qquad\sigma(I)+\sigma(J)\geq 2^{r}-1,$$ 
with equality in any of the two if and only
if there is some $1\leq l\leq r$ such that
$J=(1,\dots,1,0,\dots,0)$ ($1$'s in the first $l=|J|$ positions)
and  $I=(0,\dots,0,1,\dots,1)$ ($1$'s in the last $r-l$ positions).
In this case we have
\begin{equation}
\la t_ld^It^J\eta,[X(D,r)]\ra=(-1)^l\la\eta,[D]\ra.
\label{formulalemma}
\end{equation}
\end{lemma}
\begin{pf}
This is completely analogous to the preceding lemma.
(Formula (\ref{formulalemma}) follows from applying recursively
(3) in Lemma \ref{propdt}.)
\end{pf}

\begin{lemma}
For any sequence
$\VVV=(V_1\subset\dots\subset V_r)$ of vector bundles and $1\leq l\leq r$,
we define $R_l(\VVV):=\rk V_l$.
Let $W\in\Veqr(X(D,r))$, and define
$(V,\VVV):=\mu(r)(W)\in\PPP(X,D,r)$.
Then $$c_1(W)=c_1(V)+\sum_{1\leq l\leq r} R_l(\VVV) t_l.$$
\label{c1Wr}
\end{lemma}
\begin{pf}
Use induction on $j$ and apply at each step Corollary \ref{c1W}.
\end{pf}

Assume that $\Omega:=(\omega,a,b)\in H^2(X(D,r))$ belongs to the
Kaehler cone of $X(D,r)$. Let $W\in\Veqr(X(D,r))$ be an equivariant vector
bundle. The $\Omega$-slope of $W$ is by definition
\begin{equation}
\slope_{\Omega}W := \frac{1}{\rk W}
\left\la c_1(W) \frac{\Omega^{n+r-1}}
{(n+r-1)!},[X(D,r)]\right\ra.
\label{formslope}
\end{equation}
(Recall that $X$ has dimension $n$, so that the dimension of $X(D,r)$ is
$n+r$.)

Fix a Kaehler class $\omega\in K(X)$, and let 
$\ul=(\lambda_1,\dots,\lambda_r)$, where
$$1>\lambda_1>\dots>\lambda_r>0$$
are real numbers (the parabolic weights).
The $(\omega,\ul)$-parabolic slope of $(V,\VVV)$ is by definition
$$\parslope_{\omega,\ul}(V,\VVV)=
\frac{1}{\rk V}
\left(\left\la c_1(V)\frac{\omega^{n-1}}{(n-1)!},[X]\right\ra
+\sum_{i=1}^r R_i(\VVV)\lambda_i
\left\la \frac{\omega^{n-1}}{(n-1)!},[D]\right\ra\right).$$

Define $\beta_1:=\lambda_1$ and, for any $1<j\leq r$, 
$\beta_j:=\lambda_j/\lambda_{j-1}$. For any $\epsilon\in\RR$, let
$$\Omega(\omega,\ul,\epsilon):=\omega+\sum_{j=1}^r
\epsilon^{2^{j-1}}(d_j-\beta_j t_j).$$

\begin{theorem}
Let $W\in\Veqr(X(D,r))$, and let $(V,\VVV):=\mu(r)(W)\in\PPP(X,D,r)$.

(1) If $\epsilon>0$ is small enough, then $\Omega(\omega,\ul,\epsilon)
\in K(X(D,r))$;

(2) $\slope_{\Omega(\omega,\ul,\epsilon)}W$
is a polynomial in $\epsilon$ of the following form:
$$\sum_{J=(j,j')\atop j+j'=n+r}
\theta_J(\epsilon)\la\omega^{j}[D]^{j'},[X]\ra
+\sum_{J=(j,j')\atop j+j'=n+r-1}
\theta'_J(\epsilon)\la c_1(V)\omega^{j}[D]^{j'},[X]\ra,$$
where both sums run over pairs of nonnegative integers,
and the $\theta_J$, $\theta'_J$ are polynomials in $\epsilon$ whose
coefficients only depend
on $X$, $\omega$, $\ur$, and $\Lambda$;

(3) We then have:
$$\slope_{\Omega(\omega,\ul,\epsilon)}(W)=
\epsilon^{2^r-1}\parslope_{\omega,\ul}(V,\VVV)+O(\epsilon^{2^r}).$$
\label{calculgrau}
\end{theorem}
\begin{pf}
Statement (1) follows from Lemma \ref{KXr}. (2) follows from applying 
repeatedly Lemma \ref{propdt}.
Let us prove statement (3). To save on typing, we will write $R_l$ instead
of $R_l(\VVV)$. Combining Lemma \ref{c1Wr} with the definition
of $\Omega(\omega,\ul,\epsilon)$ we can write and develope
\begin{align*}
\slope_{\Omega(\omega,\ul,\epsilon)}W &=
\frac{1}{\rk W}\left\la
\left(c_1(V)+\sum_{1\leq l\leq r}R_lt_l\right)
\frac{\left(\omega+\sum_{j=1}^r \epsilon^{2^{j-1}}(d_j-\beta_jt_j)\right)^{n+r-1}}
{(n+r-1)!},[X(D,r)]\right\ra \\
&=\frac{1}{\rk W}\sum_{k\in\ZZ_{\geq 0},I,J\in\ZZ_{\geq 0}^r}
\theta_{k,I,J}\left
\la c_1(V)\omega^k d^I(-\beta t)^J\epsilon^{\sigma(I)+\sigma(J)},
[X(D,r)]\right\ra +\\
&+\frac{1}{\rk W}\sum_{1\leq l\leq r}
\sum_{k\in\ZZ_{\geq 0},I,J\in\ZZ_{\geq 0}^r}
\theta_{l,k,I,J}
R_l \left\la t_l\omega^k d^I(-\beta t)^J\epsilon^{\sigma(I)+\sigma(J)},
[X(D,r)]\right\ra,
\end{align*}
where the $\theta_{k,I,J}$ and $\theta_{l,k,I,J}$ are real numbers
(note that $k,I,J$ satisfy the relation $k+|I|+|J|=n+r-1$).
We now reduce mod $\epsilon^{2^r}$. By Lemmae \ref{propIJ} and \ref{propIJ2}
the only terms which can possibly be nonzero are:
\begin{enumerate}
\item in the first summation, the one with $J=0$ and $I=(1,\dots,1)$, so $k=n-1$; 
then $\theta_{k,I,J}=\left({n+r-1}\atop{r}\right)\frac{r!}{(n+r-1)!}
=\frac{1}{(n-1)!}$;
\item in the second summation,
the ones in which $J=(1,\dots,1,0,\dots,0)$ (where $1$ appears $l=|J|$ times, 
$1\leq l\leq r$) and $I=(0,\dots,0,1,\dots,1)$ (here $0$ appears $l$ times),
so again $k=n-1$;
then $\theta_{l,k,I,J}=\left({n+r-1}\atop{r}\right)
\left(r\atop i\right)\frac{i!(r-i)!}{(n+r-1)!}
=\frac{1}{(n-1)!}$.
\end{enumerate}
Taking this into account we now can write
\begin{align*}
\slope_{\Omega(\omega,\ul,\epsilon)}W &\equiv
\frac{\epsilon^{2^r-1}}{\rk W}
\left(\left\la c_1(V)\prod_{1\leq j\leq r} d_j\frac{\omega^{n-1}}{(n-1)!},
[X(D,r)]\right\ra+\right. \\
&\left.+\sum_{i=1}^r R_i\left\la t_i \prod_{1\leq j\leq i}(-\beta_jt_j)
\prod_{i<j\leq r}d_j\frac{\omega^{n-1}}{(n-1)!},[X(D,r)]\right\ra\right) 
\qquad\mod\epsilon^{2^r}.
\end{align*}
Finally, plugging in formulae (\ref{formulalemma2}) and (\ref{formulalemma}) and
taking into account that, for any $i$, $\prod_{1\leq j\leq i}\beta_j=\lambda_i$,
we deduce that
$$\slope_{\Omega(\omega,\ul,\epsilon)}W \equiv
\epsilon^{2^r-1}\parslope_{\omega,\ul}(V,\VVV)
\qquad\mod\epsilon^{2^r}.$$
This is what we wanted to prove.
\end{pf}

\section{Parabolic structures over a normal crossing divisor}
\label{ncd}
\subsection{The categories}
Let $X$ be a manifold, and let $D\subset$ be
a divisor with normal crossings.
Assume that the irreducible components $D_1,\dots,D_s$ of $D$ are smooth. 
Let us fix an $s$-tuple of nonzero natural numbers $\ur=(r_1,\dots,r_s)$.
Let $\PPP(X,D,\ur)$ be the category defined as 
follows:
\begin{enumerate}
\item The objects of $\PPP(X,D,\ur)$ are
sequences $(V,\VVV_1,\dots,\VVV_s)$, where
$V$ is a vector bundle over $X$ and where, for any $1\leq i\leq s$, 
$\VVV_i$ denotes an increasing filtration of $V|_{D_i}$ of length $r_i$:
$$\VVV_i=(0\subset V_{i,1}\subset\dots\subset V_{i,r_i}\subset V|_{D_i}).$$
\item The morphisms between two objects
$(V,\VVV_1,\dots,\VVV_s)$ and $(V',\VVV'_1,\dots,\VVV'_s)$
are the morphisms of vector bundles $\phi:V\to V'$ such that, for any
$1\leq i\leq s$, the restriction $\phi|_{D_i}$ is compatible with
the filtrations $\VVV_i$ and $\VVV'_i$,
i.e., for any $1\leq j\leq r_i$, $\phi|_{D_i}(V_{i,j})\subset V'_{i,j}$.
\end{enumerate}

We will say that a parabolic bundle $(V',\VVV'_1,\dots,\VVV'_s)$
is a parabolic subbundle of \\ $(V,\VVV_1,\dots,\VVV_s)$ if and only if
$V'\subset V$ is a subbundle and the inclusion map $\iota:V'\to V$
is a morphism between $(V',\VVV'_1,\dots,\VVV'_s)$ and $(V,\VVV_1,\dots,\VVV_s)$
in the category $\PPP(X,D,\ur)$.
In this case we will write
$$(V',\VVV'_1,\dots,\VVV'_s)\subset(V,\VVV_1,\dots,\VVV_s).$$

By Theorem \ref{mainthmr} for any $1\leq i\leq s$ there is a manifold
$\pi_i:X(D_i,r_i)\to X$ acted on by $G_i:=(\Gm)^{r_i}$ and a section
$\sigma_i:X\to X(D_i,r_i)$ of $\pi_i$.
Let us define now $X(D,\ur)$ to be the fibred product
$$X(D,\ur)=X(D_1,r_1)\times_X X(D_2,r_2)\times_X\dots\times_X
X(D_s,r_s).$$
This is smooth because the divisors $D_i$ intersect transversely.

Denote by $\Pi:X(D,\ur)\to X$ the projection. The sections $\sigma_i$
induce maps $s_i:X(D_i,r_i)\to X(D,\ur)$. If we identify 
$X(D_i,r_i)$ with $X\times_X\dots\times_X X(D_i,r_i)\times_X\dots\times_X X$
then $s_i=(\sigma_1,\dots,\sigma_{i-1},\Id,\sigma_{i+1},\dots,\sigma_s)$.
Using the sections $\sigma_i$ we also get a section $\Sigma:X\to X(D,\ur)$
of the projection $\Pi$. It is easy to check that for any $i$ we have
\begin{equation}
\Sigma=s_i\sigma_i.
\label{seccions}
\end{equation}

Consider the diagonal action of $\Gamma=G_1\times\dots\times G_s$ 
on $X(D,\ur)$. We define the category $\VeqG(X(D,\ur))$ as follows:
\begin{enumerate}
\item The objects of $\VeqG(X(D,\ur))$ are $\Gamma$-equivariant
vector bundles $W\to X(D,\ur)$ such that for any $i$ we have
$s_i^*W\in\Veqi(X(D_i,r_i))$.
\item The morphisms between two objects $W,W'$ are the $\Gamma$-equivariant
morphisms of vector bundles.
\end{enumerate}
Recall that the condition $s_i^*W\in\Veqi(X(D_i,r_i))$ translates
into a certain restriction on the weights of the action of $\Gamma$
(and hence is purely topological).

We now define a functor $M:\VeqG(X(D,\ur))\to\PPP(X,D,\ur)$. To define
the action of $M$ on objects, observe that if $W\in\VeqG(X(D,\ur))$ then
$\mu(r_i)s_i^*W=:(V_i,\VVV_i)\in\PPP(X,D_i,r_i)$ is a parabolic bundle
on $(X,D_i)$. The bundle $V_i$ is by construction the restriction
$s_i^*W|_{\sigma_i(X)}$ but by (\ref{seccions}) this is equal to
$W|_{\Sigma(X)}$. So all the bundles $V_1,\dots,V_s$ can be canonically
identified with $V:=W|_{\Sigma(X)}$ and hence the filtrations
$\VVV_1,\dots,\VVV_s$ are parabolic structures on $V$ over the divisor
$D$. We set
$$M(W):=(V,\VVV_1,\dots,\VVV_s).$$
Finally, the action of $M$ on morphisms is given by restriction on $\Sigma(X)$.

\begin{theorem} 
\label{mainthmur}
For any complex manifold $X$, any normal crossing divisor
$D\subset X$ whose irreducible components $D_1,\dots,D_s$ are smooth, 
and any sequence of nonzero natural numbers $\ur=(r_1,\dots,r_s)$, there exists 
\begin{enumerate}
\item a manifold $X(D,\ur)$ acted on by $\Gamma=(\Gm)^{|\ur|}$, 
\item an invariant projection $\Pi:X(D,\ur)\to X$ with a section
$\Sigma:X\to X(D,r)$, 
\item a full subcategory $\VeqG(X(D,\ur))$ of the category
of $\Gamma$-equivariant vector bundles on $X(D,\ur)$, and 
\item a functor $M:\VeqG(X(D,\ur))\to\PPP(X,D,\ur)$ 
\end{enumerate}
satisfying the following properties:

(A) Let $f:Y\to X$ be a map which is transverse to $D$, so that
$f^{-1}D\subset Y$ is a normal crossing divisor. Then there is an induced
map $f_{D,\ur}:Y(f^{-1}D,\ur)\to X(D,\ur)$ so that the following two diagrams
commute:
$$\xymatrix{Y(f^{-1}D,\ur)\ar[r]^-{f_{D,\ur}}\ar[d] & X(D,\ur)\ar[d] \\
Y \ar[r]^{f} & X,}\qquad\qquad
\xymatrix{\VeqG(X(D,\ur)) \ar[r]^-{f_{D,\ur}^*}\ar[d]_{M} & 
\VeqG(Y(f^{-1}D,\ur)) \ar[d]^{M} \\
\PPP(X,D,\ur) \ar[r]^{f^*} & \PPP(Y,f^{-1}D,\ur);}$$

(B) the functor $M:\VeqG(X(D,\ur))\to\PPP(X,D,\ur)$ induces an 
equivalence of categories.
\end{theorem}
\begin{pf}
(A) follows from applying inductively the commutativity of diagram
\ref{commutamu} and (B) is a consequence of Theorem \ref{mainthmr}. 
\end{pf}

\subsection{Computing the slope}

Let us assume that $X$ is Kaehler and compact.
For any $1\leq i\leq s$, let
$d(i)_1,\dots,d(i)_{r_i},t(i)_1,\dots,t(i)_{r_i}$
be the cohomology classes in $H^2(X(D_i,r_i))$ given by Lemma \ref{KXr}.
Let also $p_i:X(D,\ur)\to X(D_i,r_i)$ be the projection induced by
$\pi_1,\dots,\pi_s$, and set, for any $1\leq j\leq r_i$,
$$d_{i,j}:=p_i^*d(i)_j\text{ and }t_{i,j}:=p_i^*t(i)_j.$$
The following lemma follows from the definition of $X(D,\ur)$ and Lemma \ref{KXr}.
\begin{lemma}
The map $\pi^*:H^2(X;\RR)\to H^2(X(D,\ur);\RR)$ is injective and,
identifying $H^(X;\RR)$ with its image by this map, we have 
\begin{equation}
H^2(X(D,\ur);\RR)=H^2(X;\RR)\oplus
\bigoplus_{1\leq i\leq s}
\bigoplus_{1\leq j\leq r_i}
\RR\la d_{i,j},t_{i,j}\ra.
\end{equation}
\end{lemma}

Let us define $\sigma=2^{r_1}+\dots+2^{r_s}-s$.
In the following two lemmae we use the notations of Section \ref{slopesmooth}
(so we use standard multiindex notation; we also denote the $r_i$-tuple 
$(d_{i,1},\dots,d_{i,r_i})$ by $d_i$).
The proofs of the lemmae are easy consequences 
of Lemmae \ref{propIJ} and \ref{propIJ2}.

\begin{lemma}
\label{propIuJu}
Let $P_1,Q_1\in\ZZ_{\geq 0}^{r_1},\dots,P_s,Q_s\in\ZZ_{\geq 0}^{r_s}$.
If for some $\eta\in H^*(X)$ we have
$\la \eta\prod_{1\leq i\leq s}d_i^{P_i}t_i^{Q_i},[X(D,\ur)]\ra\neq 0$
then
($\alpha$)  for any $1\leq i\leq s$ we have 
$$|P_i|+|Q_i|\geq r_i\qquad\text{and}\qquad
\sigma(P_i)+\sigma(Q_i)\geq 2^{r_i}-1,$$
with equality in any of the two if and only if
$P_i=(1,\dots,1)$ and $Q_i=0$;

($\beta$) if $\sum_{1\leq i\leq s}\sigma(P_i)+\sigma(Q_i)\leq\sigma$
then $\sum_{1\leq i\leq s}\sigma(P_i)+\sigma(Q_i)=\sigma$,
and
$$\left\la \eta\prod_{1\leq i\leq s}d_i^{P_i}t_i^{Q_i},[X(D,\ur)]\right\ra
=\left\la \eta\prod_{1\leq i\leq s}d_i^{P_i},[X(D,\ur)]\right\ra=\la\eta,[X]\ra.$$
\end{lemma}

\begin{lemma}
\label{propIuJu2}
Let $P_1,Q_1\in\ZZ_{\geq 0}^{r_1},\dots,P_s,Q_s\in\ZZ_{\geq 0}^{r_s}$.
If for some $\eta\in H^*(X)$, $1\leq i\leq s$, $1\leq j\leq r_i$ we have
$\la t_{i,j}\eta\prod_{1\leq u\leq s}d_u^{P_u}t_u^{Q_u},[X(D,\ur)]\ra\neq 0$
then

($\alpha'$) for any $1\leq u\leq s$ we have 
$$|P_u|+|Q_u|\geq r_u\qquad\text{and}\qquad
\sigma(P_u)+\sigma(Q_u)\geq 2^{r_u}-1,$$
with equality in any of the two if and only if
\begin{enumerate}
\item either $u\neq v$, $P_u=(1,\dots,1)$ and $Q_u=0$;
\item or $u=v$, $P_u=(0,\dots,0,1,\dots,1)$ (zeroes in the first $l$ positions,
ones everywhere else) and $Q_u=(1,\dots,1,0,\dots,0)$ (ones in the first
$l$ positions, zeroes everywhere else.
\end{enumerate}

($\beta'$) if $\sum_{1\leq u\leq s}\sigma(P_u)+\sigma(Q_u)\leq\sigma$
then $\sum_{1\leq u\leq s}\sigma(P_u)+\sigma(Q_u)=\sigma$ and
$$\left\la t_{i,j}\eta\prod_{1\leq u\leq s}d_u^{P_u}t_u^{Q_u},[X(D,\ur)]\right\ra
=(-1)^l\la\eta,[X]\ra.$$
\end{lemma}

\begin{lemma}
Let $W\in\VeqG(X(D,\ur)$, and define
$(V,\VVV_1,\dots,\VVV_s):=M(W)\in\PPP(X,D,\ur)$.
Then $$c_1(W)=c_1(V)+\sum_{1\leq u\leq s}
\sum_{1\leq l\leq r_s} R_l(\VVV_u)t(u)_l.$$
\label{c1Ws}
\end{lemma}
\begin{pf}
Same idea as in Lemma \ref{c1Wr} (see loc. cit. for the definition
of $R_j$).
\end{pf}

Fix a Kaehler class $\omega\in K(X)$ and
take, for any $1\leq u\leq s$, a sequence $\ul_u=(\lambda_{u,1},\dots,
\lambda_{u,r_u})$ of real numbers satisfying
\begin{equation}
1>\lambda_{u,1}>\lambda_{u,2}>\dots>\lambda_{u,r_u}>0.
\label{condpesos}
\end{equation}
Let us write $\Lambda=(\ul_1,\dots,\ul_s)$.
The $(\omega,\Lambda)$-slope of a parabolic bundle
$(V,\VVV_1,\dots,\VVV_s)\in\PPP(X,D,\ur)$
is by definition
\begin{align*}
\parslope_{\omega,\Lambda}(V,\VVV_1,\dots,\VVV_r) &= 
\frac{1}{\rk V}
\left(\left\la c_1(V)\frac{\omega^{n-1}}{(n-1)!},[X]\right\ra + \right.\\
&+\left. \sum_{u=1}^s \sum_{i=1}^{r_u} R_i(\VVV_u)\lambda_{u,i}
\left\la \frac{\omega^{n-1}}{(n-1)!},[D_u]\right\ra\right).
\end{align*}
For any $1\leq u\leq s$, define $\beta_{u,1}=1$ and, if $1<j\leq r_u$,
$\beta_{u,j}:=\lambda_{u,j}/\lambda_{u,j-1}$.
For any $\epsilon\in\RR$, let
\begin{equation}
\Omega(\omega,\Lambda,\epsilon):=\omega+\sum_{u=1}^s
\sum_{j=1}^{r_u}\epsilon^{2^{j-1}}(d(u)_j-\beta_{u,j}t(u)_j).
\label{defOmega}
\end{equation}

\begin{theorem}
Let $W\in\VeqG(X(D,\ur))$, and let $(V,\VVV_1,\dots,\VVV_s):=M(W)
\in\PPP(X,D,\ur)$.

(1) If $\epsilon>0$ is small enough, then 
$\Omega(\omega,\Lambda,\epsilon)\in K(X(D,\ur))$;

(2) $\slope_{\Omega(\omega,\Lambda,\epsilon)}W$
is a polynomial in $\epsilon$ of the following form:
$$\sum_{Q=(j,j_1,\dots,j_s)\atop |Q|=n}
\theta_Q(\epsilon)\la\omega^{j}[D_1]^{j_1}\dots[D_s]^{j_s},[X]\ra
+\sum_{Q=(j,j_1,\dots,j_s)\atop |Q|=n-1}
\theta'_Q(\epsilon)\la c_1(V)\omega^{j}[D_1]^{j_1}\dots[D_s]^{j_s},[X]\ra,$$
where both sums run over $s+1$-tuples of nonnegative integers,
and the $\theta_Q$, $\theta'_Q$ are polynomials in $\epsilon$ whose
coefficients only depend
on $X$, $\omega$, $\ur$, and $\Lambda$;

(3) we have:
$$\slope_{\Omega(\omega,\Lambda,\epsilon)}(W)=
\epsilon^{\sigma}\parslope_{\omega,\Lambda}(V,\VVV_1,\dots,\VVV_s)
+O(\epsilon^{\sigma+1}).$$
\label{calculslope}
\end{theorem}
\begin{pf}
(1) follows from Lemma \ref{KXr}. The proof of the remaining points 
is the same as the proof of Theorem \ref{calculgrau}, but using
Lemmae \ref{propIuJu} and \ref{propIuJu} instead of Lemmae
\ref{propIJ} and \ref{propIJ2}.
\end{pf}

\begin{remark}
It should be remarked that, when chosing the parabolic weights over
an irreducible component of the divisor, one is usually allowed to 
set the smallest one equal to $0$. 
Although we do not include this possibility in our
construction (see (\ref{condpesos})), it is easy to implement it.
For example, if $\lambda_{u,r_u}=0$ for some $1\leq u\leq s$, then
we set $\beta_{u,r_u}=\epsilon$ and define $\Omega(\omega,\Lambda,\epsilon)$
as in (\ref{defOmega}). One can then check that Theorem \ref{calculslope}
remains valid.
\end{remark}

\section{Stability}
\label{sec:stab}

In this section we will assume that the Kaehler class $\omega\in K(X)$
is rational, i.e., there exists some $Q\in\NN$ so that
$\omega\in H^2(X;\ZZ[Q^{-1}])$.
This implies that for any vector bundle $V\to X$ we have
\begin{equation}
\deg V\in\ZZ[(\rk V(n-1)!Q)^{-1}].
\label{graulimitat}
\end{equation}

In the sequel, whenever we talk about a subbundle $F'$ of a vector bundle
$F\to Z$, we will implicitly mean that $V'$ is only defined over some submanifold 
$Z'\subset Z$, where $Z\setminus Z'$ is a subvariety of $Z$ of
codimension $\geq 2$ (so, strictly speaking, $V'$ is be a subbundle
of $V|_{Z'}$). This is equivalent to saying that $F'$ is reflexive subsheaf of
the sheaf of local sections of $F$. 
We remark that the degree of such subbundles is well defined (thanks to the
restriction on the codimension of $Z\setminus Z'$).
Observe also that if $W\in\VeqG(X(D,\ur))$ and 
$(V,\VVV_1,\dots,\VVV_s)=M(W)$ then, even using this extended notion of
subbundle, the functor $M$ gives a bijection between
the equivariant subbundles of $W$ and the parabolic subbundles
of $(V,\VVV_1,\dots,\VVV_s)$. This follows from the commutativity of
the diagrams in Theorem \ref{mainthmur} in the case $Y=X'$ and
$f:X'\to X$ the inclusion.

Let us recall the notions of Mumford--Takemoto (or slope) 
(semi)stability for equivariant and parabolic vector bundles.
We use the following standard notation: whenever we write a 
sentence with the word (semi)stable and the symbols $(\leq)<$ we will
mean two sentences, one with the word {\it semistable} and
the symbol $\leq$, and the other with {\it stable} and $<$.

\begin{definition}
An equivariant vector bundle $W\in\VeqG(X(D,\ur))$ is said to be
$\Omega(\omega,\Lambda,\epsilon)$-(semi)stable if and only if for 
any $\Gamma(s)$-equivariant subbundle $W'\subset W$ we have
$$\slope_{\Omega(\omega,\Lambda,\epsilon)}W'(\leq)<
\slope_{\Omega(\omega,\Lambda,\epsilon)}W.$$
\end{definition}

We remark that if $W$ is (semi)stable as a $\Gamma(s)$-equivariant
vector bundle then it is (semi)stable as a vector bundle (i.e., the
inequality between slopes holds for any subbundle of $W$, and not
only for the equivariant ones). This follows from the existence
and unicity (so in particular $\Gamma(s)$-invariance) of the
Harder--Narasimhan filtration (see \cite{GP1}).

\begin{definition}
A parabolic bundle $(V,\VVV_1,\dots,\VVV_s)\in\PPP(X,D,\ur)$ is said
to be $(\omega,\Lambda)$-(semi)stable if and only if for any parabolic
subbundle $(V',\VVV'_1,\dots,\VVV'_s)\subset(V,\VVV_1,\dots,\VVV_s)$
we have
$$\parslope_{\omega,\Lambda}(V',\VVV'_1,\dots,\VVV'_s)(\leq)<
\parslope_{\omega,\Lambda}(V,\VVV_1,\dots,\VVV_s).$$
\end{definition}

Let us fix cohomology classes $c_1\in H^2(X;\ZZ)$ and $c_2\in H^4(X;\ZZ)$.
We will call basic data the tuple $(X,D,\ur,c_1,c_2,\omega,\Lambda)$.

\begin{theorem}
There exists some $\epsilon_0>0$, depending only on the basic data,
with the following property.
Let $W\in\VeqG(X(D,\ur))$ and define $(V,\VVV_1,\dots,\VVV_s)=M(W)$.
Assume that $c_1(V)=c_1$ and $c_2(V)=c_2$. Then:

(1) If  $(V,\VVV_1,\dots,\VVV_s)$ is $(\omega,\Lambda)$-stable, 
then, for any $0<\epsilon<\epsilon_0$, $W$ is 
$\Omega(\omega,\Lambda,\epsilon)$-stable.

(2) If, for some $0<\epsilon<\epsilon_0$, 
$W$ is $\Omega(\omega,\Lambda,\epsilon)$-semistable then 
$(V,\VVV_1,\dots,\VVV_s)$ is semistable.
\label{estabilitat}
\end{theorem}

We will prove the theorem using statement (3) in Theorem \ref{calculslope}.
For that we will need to bound, uniformly over all the parabolic
subbundles of $(V,\VVV_1,\dots,\VVV_s)$, the error term.
This will be done in the next subsection, and the proof of the
theorem will be given in Subsection \ref{demostabilitat}.

\subsection{Bounding the error}

\subsubsection{}
Let $E$ be a real $2n$-dimensional vector space, 
and let $J\in\End(E)$ satisfy $J^2=-1$.
Take on $E^*$ the complex structure $-J^*\in\End(E^*)$. 
Let $\EC^*=E^*\otimes_{\RR}\CC$.
Let $a_1,\dots,a_n$ be a base of $(\EC^*){1,0}$ and let
$\ov{a}_j:=\ov{a_j}\in(\EC^*)^{0,1}$. Take on $E$ an Euclidean metric
which induces a Hermitian metric on $\EC^*$ for which
$a_1,\dots,a_n,\ov{a}_1,\dots,\ov{a}_n$ are orthogonal and satisfy
$$|a_j|^2=|\ov{a}_j|^2=2.$$
Define $\eta=\frac{\imag}{2}\sum_{i=1}^n a_i\wedge \ov{a}_i\in\Lambda^{1,1}\EC^*$.
Note that the restriction of $\eta$ to $E\subset E\otimes_{\RR}\CC$ is a real form,
and so the same thing happens to any power of $\eta$.
On the other hand, the restriction of $\eta^n/n!\in\Lambda^{n,n}\EC^*$  
to $E$ coincides with the volume form induced by the 
chosen Euclidean metric, and we have $|\eta^n/n!|=1$.

\begin{lemma}
There exists a real number $\delta_0>0$ and a constant $C>0$ so that, for any
$\theta\in\Lambda^{n-1,n-1}\EC^*$ which restricts to a real
$n-2$-form on $E$ and which satisfies $|\eta^{n-1}-\theta|<\delta_0$, we have:

(1) for any $b\in(\EC^*)^{0,1}$, 
$-\imag b\wedge \ov{b}\wedge\theta=\beta \eta^n/n!,$ where $\beta$ is a positive number;

(2) for any $g\in\Lambda^2 E^*$ we have
$g\wedge\theta=\gamma\eta^n/n!,$ where $\gamma$ is a real number satisfying
$$|\gamma|\leq C|g|.$$
\label{cotalocal}
\end{lemma}
\begin{pf}
It is clear that the map 
$\xi:(\EC^*)^{0,1}\ni b\mapsto 
-\la \imag b\wedge \ov{b}\wedge\theta,\eta^n/n!\ra$
is quadratic and takes real values. On the other hand, when
$\theta=\eta^{n-1}$, the map $\xi$ is positive definite.
Since this property is preserved by slight perturbations, (1) follows.
(2) is obvious.
\end{pf}

\subsubsection{}
\label{lemestecnics}
Let us fix from now on an Euclidean metric on 
$H^{n-1,n-1}(X;\CC)\cap H^{2n-2}(X;\RR)$, and let us denote by 
$B(\omega^{n-1},r)$ the ball
in $H^{n-1,n-1}(X;\CC)\cap H^{2n-2}(X;\RR)$ centered at $\omega^{n-1}$ 
and with radius $r$.

Recall that given a metric $h$ on a vector bundle $V$, there
exists a unique connection $A_h$ (the so called Chern connection) 
which is both compatible with $h$ and with the complex structure of $V$
(see \cite{GH}). We will denote by $F_h$ the curvature of $A_h$, and
by $\|F_h\|_{L^2}$ the $L^2$ norm of $F_h$ with respect to the metric $h$.

\begin{lemma}
Let $C_0>0$ be any real number.
There exist constants $C>0$ and $\delta>0$, depending only on $X$, $\omega$ 
and $C_0$, such that:
for any vector bundle $V\to X$ admitting a metric $h$ with
$\|F_h\|_{L^2}\leq C_0$, any subbundle $V'\subset V$, and any 
$\Theta\in B(\omega^{n-1},\delta)$, we have
$$\la c_1(V')\cap \Theta,[X]\ra\leq C.$$
\end{lemma}
\begin{pf}
Let us begin by making some observations and definitions.
There exists some $\delta_1>0$ and a map 
$\phi:B(\omega^{n-1},\delta_1)\to\Omega^{2n-2}(X)$ such that,
for any $\Theta\in B(\omega^{n-1},\delta_1)$, 
$\phi(\Theta)$ is a $2n-2$-form representing
$\Theta$.
Let $0<\delta<\delta_1$ be small enough so that for any $\Theta
\in B(\omega^{n-1},\delta)$
we have $|\phi(\Theta)-\phi(\omega^{n-1})|_{C^0}<\delta_0$, where
$\delta_0$ is the number given by Lemma \ref{cotalocal}.
Finally, let $C_0=\sup_{\Theta\in B(\omega^{n-1},\delta)}|\phi(\Theta)|_{C^0}$.

Let $V\to X$ be a vector bundle with a metric $h$ so that
$\|F_h\|_{L^2}<C_0$.
Let $V'\subset V$ be a subbundle. Using the metric $h$ we can give a
$C^{\infty}$ isomorphism $V\simeq V'\oplus V''$, where $V''=V/V'$.
By means of this splitting the $\ov{\partial}$ operator of
$V$ is the following:
$$\ov{\partial}_V=\left(\begin{array}{cc}
\ov{\partial}_{V'} & \beta \\ 0 & \ov{\partial}_{V''}\end{array}\right),$$
where $\beta\in\Omega^{0,1}(X;{V''}^*\otimes V')$ represents the element
in $H^{0,1}(X;{V''}^*\otimes V')$ corresponding to the extension
$$0\to V'\to V\to V''\to 0.$$
A standard computation gives the following formula for the curvature $F_h$
in terms of the splitting:
$$F_h=\left(\begin{array}{cc}
F_{h'}-\beta\wedge\beta^* & * \\ * & F_{h''}+\beta^*\wedge\beta\end{array}
\right),$$
where $F_{h'}$ (resp. $F_{h''}$)
$h'$ (resp. $h''$) is the restriction of $h$ to $V'$ (resp. $V''$).
Let us denote by $F_h|_{V'}\in\Omega^2(X;\ulie(V'))$ the upper left block
in the matrix. 
By the definition of $C'$ we have a pointwise bound
\begin{equation}
\Tr F_h|_{V'}\leq \rk V'|F_h|\leq \rk V |F_h|.
\label{boundcurv}
\end{equation}

Let us take some $\Theta\in B(\omega^{n-1},\delta)$, and let
$\sigma=\phi(\Theta)$  be the corresponding $2n-2$-form.
We then have 
\begin{align*}
\la c_1(V')\cap \Theta,[X]\ra
&=\int_X \frac{\imag}{2\pi}\Tr F_{h'}\wedge\sigma \\
&=\int_X \frac{\imag}{2\pi}\Tr F_{h}|_{V'}\wedge\sigma
+\int_X \frac{\imag}{2\pi}\beta\wedge\beta^*\wedge\sigma
\leq \Vol(X)\rk V \|F_h\|_{L^2},
\end{align*}
by Cauchy--Schwarz (we have used that 
$\int_X \frac{\imag}{2\pi}\beta\wedge\beta^*\wedge\sigma<0$, which follows
from (2) in Lemma \ref{cotalocal}).
This finishes the proof.
\end{pf}

\begin{lemma}
There exists a constant $C'>0$, depending only on the basic data, such that
on any $(\omega,\Lambda)$-stable parabolic bundle 
$(V,\VVV_1,\dots,\VVV_s)\in\PPP(X,D,\ur)$, with $c_1(V)=c_1$ and $c_2(V)=c_2$,
there is a metric $h$ with $\|F_h\|_{L^2}<C'$.
\end{lemma}
\begin{pf}
Let $V\to X$ be a vector bundle.
Denote by $\Lambda:\Omega^2(X)\to\Omega^0(X)$ the adjoint of the map
$\cdot\wedge\omega:\Omega^0(X)\to\Omega^2(X)$. For any metric $h$
on $V$ we have
$$\|F_h\|_{L^2}^2=\|\Lambda F_h\|_{L^2}^2-8\pi^2 ch_2(V),$$
where $ch_2(V)=\frac{1}{2}c_1(V)^2-c_2(V)$ (this follows from a simple
computation). Hence, it suffices to find a metric $h$ on $V$ with
$\|\Lambda F_h\|_{L^2}<C'',$
where $C''$ depends only on the basic data. 

Suppose that $(V,\VVV_1,\dots,\VVV_s)$ is $(\omega,\Lambda)$-stable.
It follows that there is a constant $K$ (depending only on the basic data) such
that for any subbundle $V'\subset V$ we have
\begin{equation}
\slope_{\omega}V'<\slope_{\omega}V+K.
\label{boundslope}
\end{equation}
Using the Harder--Narasimhan and Jordan--H\"older theorems (see \cite{K}) 
we get a filtration
$$V_1\subset V_2\subset\dots\subset V_j=V,$$
where each quotient $V_{i+1}/V_i$ is stable and 
$\slope(V_{i+1}/V_i)$ is bounded by a function of the basic data, thanks
to (\ref{boundslope}). 
By the Hitchin--Kobayashi correspondence (see \cite{UY})
each quotient $V_{i+1}/V_i$ admits a metric $h_{i}$ with
$$\Lambda F_{h_i}=\slope(V_{i+1}/V_i).$$
So to get the result it suffices to prove this fact: if 
\begin{equation}
0\to F'\to F\to F''\to 0
\label{extF}
\end{equation}
is an exact sequence of vector bundles
and $F'$ (resp. $F''$) admits a metric $h'$ (resp. $h''$) such that
$\|\Lambda F_{h'}\|<K'$ and $\|\Lambda F_{h''}\|<K''$ then $F$ admits
a metric $h$ satisfying $\|\Lambda F_h\|<K'+K''$. 
Now, just as in the preceding lemma, we know that $F$ is isomorphic
to the vector bundle $F'\oplus F''$ endowed with the $\ov{\partial}$-operator
$$\ov{\partial}_{F',F'',\beta}=\left(\begin{array}{cc}
\ov{\partial}_{F'} & \beta \\ 0 & \ov{\partial}_{F''}\end{array}\right),$$
where $\beta\in\Omega^{0,1}(X;{F''}^*\otimes F')$ represents the element
in $H^{0,1}(X;{F''}^*\otimes F')$ corresponding to the extension \ref{extF}.
The curvature of the Chern connection w.r.t. the metric $h'\oplus h''$ and
the above $\ov{\partial}$-operator depends continuously on $\beta$, and
tends to $F_{h'}\oplus F_{h''}$ as the norm of $\beta$ goes to zero.
But if we substitute $\beta$ by $\lambda\beta$ for any $\lambda\in\CC$
we do not change the isomorphism class of the bundle. So it suffices
to take $\lambda$ small enough and we are done.
\end{pf}

\begin{corollary}
\label{bounderror}
There exists some constant $C''>0$, depending only on
the basic data such that, for any $(\omega,\Lambda)$-stable
$(V,\VVV_1,\dots,\VVV_s)$ satisfying $c_1(V)=c_1$ and $c_2(V)=c_2$, 
any subbundle $V'\subset V$, and any
$s+1$-tuple of nonnegative integers $J=(j,j_1,\dots,j_s)$ such that
$|J|=n-1$, we have
$$|\la c_1(V')\omega^j[D_1]^{j_1}\dots[D_s]^{j_s},[X]\ra|\leq
C''(1+|\la c_1(V')\omega^{n-1},[X]\ra|).$$
\end{corollary}
\begin{pf}
Let $\epsilon>0$ be small enough so that for any $J=(j,j_1,\dots,j_s)$
we have 
$\omega^{n-1}\pm\epsilon\omega^j[D_1]^{j_1}\dots[D_s]^{j_s}
\in B(\omega^{n-1},\delta)$.
The preceding two lemmae imply that for any $J=(j,j_1,\dots,j_s)$
$$\pm\epsilon\la c_1(V')\omega^j[D_1]^{j_1}\dots[D_s]^{j_s},[X]\ra
\leq C-\la c_1(V')\omega^{n-1},[X]\ra,$$
from which the result follows easily.
\end{pf}

\subsection{Proof of Theorem \ref{estabilitat}}
\label{demostabilitat}
In all this subsection, whenever we say {\it for small enough $\epsilon$},
we will implicitly mean {\it depending only on the basic data}.

(1) Let us take an equivariant vector bundle $W\in\VeqG(X(D,\ur))$,
and define $(V,\VVV_1,\dots,\VVV_s)=M(W)$. 
Assume that $(V,\VVV_1,\dots,\VVV_s)$ is $(\omega,\Lambda)$-stable.
Let $R=\rk V$. By (\ref{graulimitat}), for any subbundle $V'\subset V$
we have $\deg V'\in\ZZ[(R!(n-1)!Q)^{-1}]$. This, together with the stability
condition implies the existence of some $\alpha>0$ such that for any
parabolic subbundle $(V',\VVV'_1,\dots,\VVV'_s)\subset(V,\VVV_1,\dots,\VVV_s)$
we have
\begin{align}
\parslope_{\omega,\Lambda}(V',\VVV'_1,\dots,\VVV'_s)
&+\alpha|\parslope_{\omega,\Lambda}(V',\VVV'_1,\dots,\VVV'_s)| \notag \\
&\leq \parslope_{\omega,\Lambda}(V,\VVV_1,\dots,\VVV_s)-\alpha.
\label{desstable}
\end{align}
We want to prove that, for small enough $\epsilon$, 
$W$ is $\Omega(\omega,\Lambda,\epsilon)$-stable. In other words:
letting $W'\subset W$ be a $\Gamma(s)$-equivariant subbundle, we
want to check that for small enough $\epsilon$
$$\slope_{\Omega(\omega,\Lambda,\epsilon)}W'<
\slope_{\Omega(\omega,\Lambda,\epsilon)}W.$$
Let $(V',\VVV'_1,\dots,\VVV'_s)=M(W')$. By Theorem \ref{mainthmur}
$(V',\VVV'_1,\dots,\VVV'_s)\subset (V,\VVV_1,\dots,\VVV_s)$, so
the inequality (\ref{desstable}) holds.
By (3) in Theorem \ref{calculslope} for small enough $\epsilon>0$ we have
$$\slope_{\Omega(\omega,\Lambda,\epsilon)}W
>\epsilon^{\sigma}\left(\parslope_{\omega,\Lambda}(V,\VVV_1,\dots,\VVV_s)-
\frac{\alpha}{2}\right).$$
On the other hand, by (2) and (3) in Theorem \ref{calculslope} we have
\begin{align*}
\slope_{\Omega(\omega,\Lambda,\epsilon)}W' &=
\epsilon^{\sigma}(\parslope_{\omega,\Lambda}(V',\VVV'_1,\dots,\VVV'_s)) \\
&+\epsilon^{\sigma+1}
\left(\sum_J\epsilon^{\deg(J)}\theta_JE_J+\sum_J\epsilon^{\deg(J)}\theta'_JF_J
\right),
\end{align*}
where $J$ denotes tuples $(j,j_1,\dots,j_s)$, $\deg(J)\geq 0$ for any $J$,
$\theta_J$ and $\theta'_J$ are numbers which only depend on the basic
data, the $E_J$'s are of the form $\la \omega^j[D_1]^{j_1}\dots[D_s]^{j_s},[X]\ra$
(so they also depend only on basic data) and where, finally, 
the $F_J$'s are of the form
$$F_J=\la c_1(V')\omega^j[D_1]^{j_1}\dots[D_s]^{j_s},[X]\ra.$$
By Corollary \ref{bounderror} we have bounds
$$|F_J|\leq C''(1+|\slope_{\omega}V'|),$$
($C''$ depends only on basic data) so, taking into account
that $$|\slope_{\omega}V'|\leq
|\parslope_{\omega,\Lambda}(V',\VVV'_1,\dots,\VVV'_s)|+C'''$$
for some $C'''$ depending on the basic data, it follows that, for
small enough $\epsilon$, we have
$$\slope_{\Omega(\omega,\Lambda,\epsilon)}W'<
\epsilon^{\sigma}
\left(\parslope_{\omega,\Lambda}(V',\VVV'_1,\dots,\VVV'_s)
+\alpha|\parslope_{\omega,\Lambda}(V',\VVV'_1,\dots,\VVV'_s)|+
\frac{\alpha}{2}\right),$$
and we are done.

The proof of (2) follows exactly the same lines, so is left to the reader.

\section{The case $\dim_{\CC}X=1$ and $r=1$} 
\label{sec:dim1}

Suppose that $X$ is a compact Riemann surface (so that $D$ is a finte 
set of points) and that $r=1$. Under these hypothesis, we obtain a 
stronger result than in the general case. 
Before stating it, let us recall that we have
$H^2(X_D)=H^2(X)\oplus\RR\la d,t\ra$  (by Lemma \ref{H2XDj})
and that $td=d^2=t p^*\alpha=0$
(here $\alpha\in H^2(X)$ is arbitrary and $p:X_D\to X$ denotes
the projection), $\la d(p^*\alpha),[X_D]\ra=\la\alpha,[X]\ra$, and
$t^2=-1$ (see Lemma \ref{propdt}).

\begin{theorem}
Let $W\in\Veq(X_D)$ and let $(V,V')=\mu(W)\in\PPP(X,D)$. 
Let $0<\alpha<1$ and let $\omega\in H^2(X)$ be a Kaehler class.
Define $\Omega=p^*\omega-\alpha t$. We then have

(1) $\Omega$ is a Kaehler class of $X_D$,

(2) $\slope_{\Omega}(W)=\parslope_{\omega,\alpha}(V,V')$,

(3) $W$ is $\Omega$-(semi)stable $\Longleftrightarrow$ $(V,V')$ is
$(\omega,\alpha)$-parabolic (semi)stable,

(4) fix some $d\in H^2(X;\ZZ)$, some nonzero $r\geq r'\in\NN$; let 
$\MMM_{d,r,r',\alpha}(X)$ be the moduli space of $(\omega,\alpha)$-stable
parabolic vector bundles $(V,V')$ over $(X,D)$ satisfying 
$c_1(V)=d$, $\rk V=r$, $\rk V'=r'$;
let $\MMM_{\Omega}(X_D)$ be the moduli space of vector bundles
$W\to X_D$ with Chern character $\ch W=p^*(r+d)+tr'$;
then $\Gm$ acts algebraically on $\MMM_{\Omega}(X_D)$ and 
$\MMM_{d,r,r',\alpha}$ can be identified with some of the connected 
components of the fixed point set of this action.
\label{millorcas}
\end{theorem}
\begin{pf}
(1) follows from Theorem \ref{KXX}; (2) can be checked by computing
both sides, using the relations between $d$ and $t$ which we recalled
above; (3) is a consequence of (2); finally, (4) follows from (2),
the ideas in Section 5.1 of \cite{GP2},
and the fact that Theorem \ref{mainthm} works also for families of
vector bundles --- this is true because the only technical tool which
was used in the proof was Riemann's extension theorem, which is
certainly true for pairs $(X\times S,D\times S)$ for any scheme
$S$ of finite type.
\end{pf}

We observe that not only Theorem \ref{mainthm} but its most general
version, Theorem \ref{mainthmur}, is valid for families of vector
bundles (by exactly the same argument as before). Unfortunately, this
does not allow to identify the moduli spaces of equivariant and
parabolic vector bundles, since the best thing we could do in general
was to identify the notion of 
$(\omega,\Lambda)$-parabolic stability to that of
$\Omega(\omega,\Lambda,\epsilon)$-stability as $\epsilon\to 0$. 
If we knew that the notion of $\Omega(\omega,\Lambda,\epsilon)$-stability 
for vector bundles over $X(D,\ur)$ does not change when $\epsilon>0$ is 
small enough, then we could identify the moduli spaces. But this does not
seem to be the case in general (see for example \cite{Q,S}).


\begin{thebibliography}{ABCD}\frenchspacing\smallbreak

\bibitem[A]{A} L. {\'A}lvarez C{\'o}nsul,
Dimensional reduction and Hitchin--Kobayashi correspondence,
{\em Ph. D. Thesis}, Universidad Aut{\'o}noma de Madrid, 
May 2000.

\bibitem[AGP]{AGP}
L. {\'A}lvarez C{\'o}nsul, O. Garc{\'\i}a--Prada, 
Dimensional reduction, $\SL(2,\CC)$-equivariant bundles and stable 
holomorphic chains, to appear in {\em Int. J. Math.}.

\bibitem [B]{B} I. Biswas, Parabolic bundles as orbifold bundles,
{\em Duke Math. J.} {\bf 88} (1997) 305--325.

\bibitem[BDGW]{BDGW}
        S.B. Bradlow, G. Daskalopoulos,  O. Garc{\'\i}a--Prada
        and R. Wentworth, Stable augmented bundles over Riemann surfaces.
        {\em Vector Bundles in Algebraic Geometry}
        (Durham 1993), eds. N.J. Hitchin, P.E. Newstead and W.M.
        Oxbury, LMS Lecture Notes Series, 208 (1995) 15--68,
        Cambridge University Press.

\bibitem[BGM]{BGM}
S.B. Bradlow, O. Garc{\'\i}a--Prada, I. Mundet i Riera,
A Hitchin--Kobayashi correspondence for subgroups of the gauge group,
in preparation.

\bibitem[CS]{CS}
J.B. Carrell, A.J. Sommese, $\CC^*$-actions, {\em Math. Scand.} {\bf
  43} (1978) 49--59; Corrections, {\em ibid.}, {\bf 53} (1983) 32.

\bibitem[GP1]{GP1} O. Garc\'{\i}a--Prada, Invariant connections and
vortices, {\em Comm. Math. Phys.} {\bf 156} (1993) 527--546.

\bibitem[GP2]{GP2} O. Garc\'{\i}a--Prada, Dimensional reduction of
stable bundles, vortices and stable pairs, {\em Int. J. of Math.}
{\bf 5} (1994) 1--52.

\bibitem[GM]{GM} I. Gelfand, Yu. Manin, Homological Algebra, Springer Verlag,
1996.

\bibitem[GH]{GH} P. Griffiths, J. Harris, Principles of Algebraic Geometry,
John Wiley \& Sons, 1978.

\bibitem[GLS]{GLS}
V. Guillemin, E. Lerman, S. Sternberg,
{\em Symplectic Fibrations and Multiplicity Diagrams},
Cambridge University Press, 1996.

\bibitem[K]{K} S. Kobayashi, Differential Geometry of Complex Vector
Bundles, Princeton Univ. Press, 1987.

\bibitem[M]{M} I. Mundet i Riera, Yang--Mills--Higgs theory for
  symplectic fibrations, {\em Ph. D. Thesis}, Universidad Aut{\'o}noma
  de Madrid, April 1999.

\bibitem[Q]{Q} Z. Qin, Equivalence classes of polarizations and moduli 
spaces of sheaves, {\em J. Diff. Geom.} {\bf 37} (1993), no. 2,
397--415. 

\bibitem[S]{S} A. Schmitt, 
Walls for Gieseker semistability and the Mumford-Thaddeus principle for moduli spaces of sheaves over higher dimensional
bases, {\em Comment. Math. Helv.} {\bf 75} (2000), no. 2, 216--231. 

\bibitem[UY]{UY} K. Uhlenbeck, S.T. Yau, On the existence of Hermitian--Yang--Mills
connections in stable vector bundles, {\em Comm. Pure. Appl. Math.}
{\bf 39-S} (1986) 257---293.

\end{thebibliography}
\end{document}